\documentclass[11pt]{article}
\usepackage{amsfonts,latexsym,amstext} 
\usepackage{amsmath} 
\usepackage{amssymb} 
\usepackage{comment} 
\usepackage{amsthm} 
\usepackage[latin1]{inputenc} 
\usepackage{color} 
 
\theoremstyle{plain} 
\newtheorem{theorem}{Theorem}[section] 
\newtheorem{lemma}[theorem]{Lemma} 
\newtheorem{proposition}[theorem]{Proposition} 
\newtheorem{corollary}[theorem]{Corollary} 

\newtheorem{definition}[theorem]{Definition} 
\newtheorem{hypothesis}[theorem]{Hypothesis} 
\theoremstyle{remark} 
\newtheorem{remark}[theorem]{Remark} 

\allowdisplaybreaks 
 
\makeatletter 
\@addtoreset{equation}{section} 
\def\theequation{\thesection.\arabic{equation}} 
\makeatother 
 
\def\qed{{\hfill\hbox{\enspace${ \square}$}} \smallskip} 
\def\sqr#1#2{{\vcenter{\vbox{\hrule height .#2pt \hbox{\vrule 
 width .#2pt height#1pt \kern#1pt \vrule 
width .#2pt} \hrule height .#2pt}}}} 
\def\square{\mathchoice\sqr54\sqr54\sqr{4.1}3\sqr{3.5}3}

\def\ds{\begin{displaystyle}} 
\def\eds{\end{displaystyle}} 
 
\def\<{\langle } 
\def\>{\rangle }

\def\R{\mathbb R} 
\def\N{\mathbb N}

\def\E{\mathbb E} 
\def\P{\mathbb P}

\newcommand{\sper}[1]{\mathbb{E} \left[ #1 \right]}                               
 
\DeclareMathAlphabet{\mathonebb}{U}{bbold}{m}{n}                           %
\newcommand{\one}{\ensuremath{\mathonebb{1}}}                               
 
\title{Special weak Dirichlet processes and 
 BSDEs driven by a random measure} 

\usepackage{authblk}

\author[1,2]{Elena Bandini\thanks{elena.bandini@polimi.it}} 
\author[2]{Francesco Russo\thanks{francesco.russo@ensta-paristech.fr}} 
\affil[1]{Politecnico di Milano, Dipartimento di Matematica, via Bonardi 9, 20133 Milano, Italy} 
\affil[2]{ENSTA ParisTech, Universit\'e Paris-Saclay, Unit\'e de Math\'ematiques appliqu\'ees, 828, Boulevard des Mar\'echaux, F-91120 Palaiseau, France}

\date{} 
 
\begin{document} 
 
\allowdisplaybreaks 
\maketitle

\begin{abstract} 
This paper considers a forward BSDE  driven by a random 
measure, when the underlying forward process $X$ is a 
special semimartingale, or even more generally,  
a special weak Dirichlet process. 
Given a solution $(Y,Z,U)$, generally $Y$ appears to be 
of the type $u(t,X_t)$ where $u$ is a deterministic function. 
In this paper we identify $Z$ and $U$ in terms of $u$  
applying  stochastic calculus with respect to weak Dirichlet processes.

\end{abstract} 

 
\section{Introduction} 
 
This paper considers a forward BSDE  driven by a random 
measure, when the underlying forward process $X$ is  a
special semimartingale, or even more generally,  
a special weak Dirichlet process. 
Given a solution $(Y,Z,U)$, often $Y$ appears to be 
of the type $v(t,X_t)$ where $v$ is a deterministic function. 
In this paper we identify $Z$ and $U$ in terms of $v$  
applying  stochastic calculus with respect to weak Dirichlet processes. 
 
Indeed the employed techniques perform  the calculus with respect to (special) 
 weak Dirichlet processes developed in \cite{BandiniRusso1}, extending the techniques established in the continuous framework  in \cite{gr}, \cite{gr1}.  
Given some filtration $(\mathcal F_t)$, we recall that a special weak Dirichlet process is a process of the type $X= M+ A$, where $M$ is an $(\mathcal F_t)$-local martingale and $A$ is an $(\mathcal F_t)$-predictable orthogonal process, see Definition 
5.6 in \cite{BandiniRusso1}. 
When $A$ has bounded variation, then $X$ is a special $(\mathcal F_t)$-semimartingale. 
That calculus has two important features: 1) the decomposition 
of a special weak Dirichlet process is unique, see Proposition \ref{Sec:WD_P_unique_decomp};
2) 
there is a chain rule (in substitution of It\^o's formula)
allowing to expand $v(t, X_t)$, where $X$ is a special weak Dirichlet process of finite quadratic variation and $v$ is of class $C^{0,1}$, fulfilling some technical assumption, see   Theorem 
\ref{T_C1_dec_specialweak_Dir}. 
	If we know a priori that $v(t,X_t)$ is the sum of a bounded variation process and a continuous $(\mathcal F_t)$-orthogonal process, then the chain rule does not require any differentiability on $v$; in that case, no assumptions are required on the c\`adl\`ag process $X$, see Proposition \ref{P_C00_chain_rule}. 
 
As we have already mentioned, we will focus on forward BSDEs, which constitute a particular case of BSDEs in their general form. 
BSDEs  have been deeply studied   since the seminal paper 
\cite{pardouxpeng}. 
In \cite{pardouxpeng}, 
as well as in many subsequent papers, the standard Brownian motion 
is the driving process (Brownian context) and  the concept of 
BSDE is based on a non-linear martingale 
representation theorem with respect to the corresponding 
Brownian filtration. 
A recent monograph on the subject is \cite{rascanu2014stochastic}. 
BSDEs driven by processes with jumps have also been  investigated: 
two classes of such equations appear in the literature. The first one 
 relates to BSDEs where the Brownian motion is replaced by 
a general c\`adl\`ag martingale $M$, see, among others, \cite{Buckdahn93}, \cite{ElKarouiHuang}, \cite{CarboneFerrSantacroce}, \cite{ceci_cretarola_russo}. 
An alternative version of BSDEs with a discontinuous driving term 
is the one associated to an integer-valued random measure  $\mu$, with corresponding compensator $\nu$. In this case the BSDE is driven by a continuous martingale $M$ and 
 a compensated random measure $\mu -\nu$. In that equation naturally appears 
a purely discontinuous martingale which is a stochastic integral with respect to 
$\mu -\nu$, see, e.g.,  \cite{xia}, \cite{BP}, \cite{TaLi}, \cite{BandiniFuhrman}, \cite{BandiniPDMPs}, \cite{BandiniConfortola}. 
A recent monograph on BSDEs driven by Poisson random measures is \cite{delong2013backward}. 
Connections between the martingale  and the random measure driven  BSDEs are 
illustrated  by  \cite{JeanblancManiaSantacroceSchweizer}.

In this paper we will focus on 
 BSDEs driven by a compensated random measure $\mu-\nu$ (we will use the one-dimensional formalism for simplicity). We will not ask $\mu$ to be  quasi-left-continuous, i.e.
$\mu(\{S \} \times \R)=0$ on $\{S < \infty \}$, for every predictable time $S$, see  the definition in
  Theorem (4.47) in \cite{jacod_book}.
Let $T>0$ be a finite horizon time. Besides $\mu$ and $\nu$ appear three driving random elements: a continuous martingale $M$, 
a non-decreasing adapted continuous process $\zeta$ and a predictable random measure 
$\lambda$ on $\Omega \times [0,T] \times \R$, equipped with the usual product $\sigma$-fields. 
Given a square integrable random variable $\xi$, and two measurable functions $\tilde g:\Omega \times [0,T]   \times \R^2 \rightarrow \R$, 
$ \tilde f: \Omega \times [0,T] \times \R^3 \rightarrow \R$, 
the equation takes the  form 
\begin{align}\label{GeneralBSDE} 
Y_t &= \xi + \int_{]t,\,T]} \tilde g(s,\,Y_{s-},\,Z_s)\,  d \zeta_s + 
\int_{]t,\,T]\times \R} \tilde f(s,\,e,\,Y_{s-},\,U_{s}(e))\,  \lambda(ds\,de)\nonumber\\ 
&\quad- \int_{]t,\,T]} Z_s \, dM_s - \int_{]t,\,T]\times \R} U_s(e)\,(\mu-\nu)(ds\,de). 
\end{align} 
As we have anticipated before, 
the unknown of \eqref{GeneralBSDE} is a triplet $(Y,Z,U)$ where $Y, Z$ are adapted 
 and $U$ is a predictable random field. 
The Brownian context of Pardoux-Peng appears as a particular case, setting 
$\mu= \lambda=0$, $\zeta_s \equiv s$. 
There, $M$ is a standard Brownian 
motion and $\xi$ is measurable with respect to the Brownian $\sigma$-field 
at terminal time. 
 In that case  the unknown can be reduced to $ (Y,Z)$, since   $U$ can be arbitrarily chosen. 
Another significant  subcase of \eqref{GeneralBSDE} arises when only the 
  purely discontinuous driving term 
appears, i.e. $M$ and 
$\zeta$ vanish.

The standard situation in the literature corresponds to the case when $\mu$ is quasi-left-continuous, see e.g. \cite{BaBuPa}, \cite{TaLi} when $\mu$ is a Poisson random measure, and \cite{Be} for a  random measure $\mu$ more general than the Poisson one, whose  compensator is absolutely continuous with respect to a deterministic measure. In the purely discontinuous subcase,
 for instance when $\mu$ is the jump measure of  a marked point process,
the well-posedness of the related BSDE can be settled
by an iterative method,
 see
 \cite{CFJ}. 
Existence and uniqueness for 
BSDEs driven by a random measure $\mu$ which is not necessarily quasi-left-continuous are very recent, 
and have been 
 discussed   in  \cite{BandiniBSDE} in the purely discontinuous case, and in  a slightly different context 
by   \cite{CohenElliott}, for  BSDEs driven by a countable sequence of square-integrable 
martingales.
 
When the random dependence 
of   $\tilde f$ and $\tilde g$  is provided by a Markov solution $X$ of a 
 forward  SDE, and $\xi$ is a real function of $X$ at the terminal time $T$, then the   BSDE \eqref{GeneralBSDE} 
 is called forward BSDE, as mentioned 
 at the beginning. 
Forward BSDEs  generally constitute  stochastic representations of a 
partial integro-differential equation (PIDE). 
In the Brownian case, 
when $X$ is the solution of a classical SDE with 
diffusion coefficient $\sigma$, then the PIDE reduces to a semilinear parabolic PDE. 
If  $v:[0,T] \times \R \times \R$ is a 
classical (smooth) solution of the mentioned  
PDE, 
then $Y_s = v(s,X_s)$, $Z_s = \sigma(s,X_s)\, \partial_x v(s,X_s)$, 
 generate  a solution 
to the forward BSDE, 
see e.g. \cite{Peng91}, \cite{PardouxPeng92}, \cite{Peng92b}. 
In the general case when 
the forward  BSDEs are 
also  driven by random measures, 
similar results have been established, for instance 
by \cite{BaBuPa}, for the jump-diffusion case, and by \cite{CoFu-m},  for the purely discontinuous case, in particular when no Brownian noise appears. 
In the context of martingale driven forward BSDEs, a first approach to the probabilistic representation 
has been carried on in \cite{LaachirRusso}. 
 
Conversely, 
solutions of forward BSDEs generate solutions  of PIDEs in the viscosity sense. 
More precisely, for each given couple $(t,x) \in [0,\,T]\times \R$, consider an  underlying process $X$ given by 
the solution $X^{t,x}$ of an SDE starting at $x$ at time $t$. 
Let $(Y^{t,x},Z^{t,x},U^{t,x})$ be a family of solutions of the forward BSDE. 
In that case, under reasonable general assumptions, the function $v(t,x):= Y_t^{t,x}$ is 
a viscosity solution of the related PIDE. 
A demanding task consists in  characterizing the couple 
$(Z,U):=(Z^{t,x},U^{t,x})$, in term of $v$; this  is generally called the \emph{identification problem} of $(Z,U)$. 
In the Brownian context, it was for instance the object of  \cite{FuhrmanTessitore}: the authors show 
that if $v \in C^{0,1}$, 
then $Z_s = \sigma(s,X_s)\,\partial_x v(s,X_s)$; under more general assumptions, 
the authors also associate $Z$ with a  generalized gradient of $v$.
At our knowledge, in the general case, the problem of the identification of the 
martingale integrands couple $(Z,U)$  has not been deeply investigated, except for particular situations, as for instance  the one   treated in  \cite{CFJ}: this problem was faced in 
%
\cite{CoFu-m}. 

A motivating PIDE is the one of Hamilton-Jacobi-Bellman related to stochastic optimal control problems, when the underlying is a general jump process. The solution to the identification problem in the related BSDE can be  useful to determine feedback control strategies in verification theorems.
Those verification theorems have the advantage of requiring less regularity of the value function than the classical ones, which need, instead, a time-space $C^{1,2}$ regularity.
Another  possible application 
concerns hedging in incomplete markets in mathematical finance. This can be treated making explicit the  so called F\"ollmer-Schweizer decomposition of the payoff related to a contingent claim 
via mean-variance hedging; that decomposition is 
strictly related to a specific BSDE, as for instance illustrated in \cite{JeanblancManiaSantacroceSchweizer}.
When the underlying is a general Markov process and the contingent claim is of vanilla type, solving  the identification problem gives us suitable techniques to discuss  mean-variance hedging.

In the present paper we discuss the over-mentioned  identification problem in a more general framework.
In particular we have formulated a set of hypotheses
which include the existing results in the literature.
In Section \ref{Sec_application_BSDE} 
 we state  Theorems \ref{P_ident} and \ref{P_ident_C0}, which are the main results of the article.
 If $Y$ together with $(Z,U)$ constitutes a solution of a BSDE and there exists a function $v$ with some minimal regularity such that $Y_t = v(t, X_t)$, those two theorems allow to solve the identification problem.
Their proof   makes  use in essential way of both features of the calculus related  to weak Dirichlet processes.
Indeed, first, 
supposing that 
  $Y$ is a suitable $C^{0,1}$-deterministic function $v$ of 
the underlying process $X$, which is a  special semimartingale $X$, 
related in a specific way to the  random measure $\mu$, 
we apply the chain rule in Theorem \ref{T_C1_dec_specialweak_Dir}. 
 In particular $Y$ will be a special  weak Dirichlet process
with prescribed local martingale part. 
Second, 
 $Y$ can be decomposed using the fact that, together with $(Z,U)$, it solves the BSDE.
So the uniqueness of decomposition of the special weak Dirichlet process $Y$ recalled in Proposition \ref{Sec:WD_P_unique_decomp}
allows us to  identify the couple $(Z,U)$. This is the object of 
 Theorem \ref{P_ident}.
On the other hand, in the purely discontinuous  framework
 we make use instead of the chain rule 
Proposition \ref{P_C00_chain_rule}. This does not even ask 
 $X$  to be a   special weak Dirichlet process, 
 provided we have some a priori information on the structure of $v(t,X_t)$.
  At this point, by a similar argument, Theorem \ref{P_ident_C0} also allows to tackle and solve the
identification problem.  

Section \ref{Sec_Appl} is devoted to concrete applications
in the following situations.
\begin{enumerate}
\item The case when the BSDE is driven by a
Wiener process and a Poisson random measure and the underlying
process $X$ is a jump-diffusion. 
\item A case when the
BSDE is driven by a quasi-left-continuous random measure 
and $X$ is a non-diffusive Markov processes. This
case has also been treated with a different technique by
\cite{CoFu-m}.
\item A particular case of a BSDE 
 driven by a non quasi-left-continuous random measure.
\end{enumerate}





The paper is organized as follows. 
In Section \ref{Sec_2.1} we  fix the  notations and we make some technical observations on  stochastic integration with respect to a random measure;  
in Section \ref{Sec_2.2} we introduce our basic set of hypotheses, and we provide some related technical results, which are proved in the Appendix.
Section \ref{Sec_application_BSDE} 
 is devoted to solve the identification  problem. As mentioned earlier, applications are provided in Section \ref{Sec_Appl}.


\section{Notations and preliminaries}\label{Sec_NOt_Prel}
In what follows, we are given a probability space $(\Omega,\mathcal{F},\P)$ 
a positive horizon $T$ 
and a filtration 
$(\mathcal{F}_t)_{t \geq 0}$, satisfying the usual conditions. Let $\mathcal F = \mathcal F_T$. 
Given a topological space $E$, in the sequel $\mathcal{B}(E)$ will denote 
the Borel $\sigma$-field associated with $E$. 
$\mathcal{P}$ (resp. $\mathcal{\tilde{P}}=\mathcal{P}\otimes \mathcal{B}(\R)$) will denote the predictable $\sigma$-field on $\Omega \times [0,T]$ (resp. on $\tilde{\Omega}= \Omega \times [0,T]\times \R$). 
Analogously, we set $\mathcal{O}$ (resp. $\tilde{\mathcal{O}} =\mathcal{O} \otimes \mathcal{B}(\R)$) as the optional   $\sigma$-field on $\Omega \times [0,T]$ (resp. on $\tilde{\Omega}$). 
Moreover, 
$\tilde {\mathcal{F}}$ will be $\sigma$-field $\mathcal{F} \otimes {\mathcal B}([0,T] \times \R)$, and 
we will  indicate by $\mathcal{F^{\P}}$ 
the completion of $\mathcal F$ with the $\P$-null sets. We set 
$\mathcal{\tilde F}^{\P} =  \mathcal{F}^{\P} \otimes {\mathcal B}([0,T] \times \R)$. 
By default, all the stochastic processes will be considered with parameter $t \in [0,\,T]$. 
By convention, any c\`adl\`ag process defined on $[0,\,T]$ is  extended to $\R_+$ by continuity.
A random set $A \subset \tilde \Omega$ is called evanescent if the set $\{\omega: \exists t \in \R_+ \textup{ with } (\omega, t) \in A\}$ is $\P$-null. Generically, all the equalities of random sets will be intended up to an evanescent set.

 A bounded variation process $X$ on $[0,\,T]$  will be said to be with integrable variation if the expectation of its total variation is finite. 
$\mathcal{A}$ (resp. $\mathcal{A}_{\textup{loc}}$) will denote  the collection of all adapted processes with   integrable variation (resp.  with locally integrable variation), and    $\mathcal{A}^+$ (resp $\mathcal{A}_{\textup{loc}}^+$)  the collection of all adapted integrable increasing (resp. adapted locally integrable)  processes. 
The significance of {\it locally} is the usual one which refers 
to  localization by stopping times, see e.g. (0.39) of 
\cite{jacod_book}. 
We will indicate by 
$C^{0,1}$ 
the space of all functions 
$$ 
u: [0,\,T]\times \R \rightarrow \R, \quad (t,x)\mapsto u(t,x) 
$$ 
that are continuous  together their derivative 
$\partial_x u$. 

\subsection{Stochastic integration with respect to integer-valued random measures} \label{Sec_2.1}
The concept of random measure  
will be extensively used 	throughout   the paper. 
For a  detailed  discussion on this topic  and the unexplained  notations see 
Chapter I and Chapter II, Section 1, in \cite{JacodBook}, Chapter III in \cite{jacod_book},  and  Chapter XI, Section 1, in \cite{chineseBook}.
In particular, if $\mu$ is a random measure on $[0,\,T]\times \R$, for any measurable real function $H$ defined on $\tilde \Omega$, one denotes $H \star \mu_t:= \int_{]0,\,t] \times \R} H(\cdot, s,x) \,\mu(\cdot, ds \,dx)$, at least when the right-hand side is strictly greater than $- \infty$.

In the sequel of the section  $\mu$ will be  an integer-valued random measure on $[0,\,T] \times \R$, and $\nu$ a "good" version of the compensator of $\mu$, as constructed in point (c) of 
Proposition 1.17, Chapter II, in \cite{JacodBook}.
Set 
$D=\{(\omega, t): \mu(\omega, \{t\}\times \R)>0\}$, 
and
$$J=\{(\omega, t): \nu(\omega, \{t\}\times \R)>0\},\quad 
K=\{(\omega, t): \nu(\omega, \{t\}\times \R)=1\}.$$
We define $\nu^d:= \nu\,\one_{J}$ and $\nu^c:= \nu\,\one_{J^c}$.
\begin{remark}\label{R_pred_supp}
	$J$ is the predictable support of $D$, see  Proposition 1.14, Chapter II,  in \cite{JacodBook}.
	The definition of  predictable support of a random set is given in  
Definition 2.32, Chapter I, in \cite{JacodBook}.
From that it follows that  
$\one_J= 
{}^p{}(\one_D)$. 
	$K$ is the largest predictable subset of $D$, 	
	 see  
	 Theorem 11.14 in \cite{chineseBook}. 
Since $K$ is predictable, 
we have 
${}^p{} 
(\one_K ) 
= 
\one_K$. 
$J$ is a thin set, see Proposition 2.34, Chapter I,  in \cite{JacodBook}.
\end{remark} 

\begin{remark}\label{R_H_mu} 
	\begin{itemize} 
		\item[(i)]$\nu$ admits a disintegration of the type 
		\begin{equation}\label{nu_dis} 
		\nu(\omega, ds\,de)=d A_s(\omega) \,\phi(\omega, s,de), 
		\end{equation} 
		where  $\phi$ is a random kernel from $(\Omega \times [0,\,T], \mathcal P)$ into $(\R, \mathcal B(\R))$ and  $A$ is a right-continuous nondecreasing predictable process, such that $A_0=0$, see 
		for instance Remark 4.4 in \cite{BandiniBSDE}. 
		\item[(ii)] 
		Given $\nu$ in the form \eqref{nu_dis}, then 
		the process $A$ is continuous if and only if, up to an evanescent set,  $D=   \cup_n [[T^i_n]]$, $[[T^i_n]] \cap [[T^i_m]] = \emptyset$, $n \neq m$, where $(T^i_n)_n$ are totally inaccessible times, see, e.g., Assumption (A) in  \cite{CFJ}. 
	\item [(iii)] We recall that a totally inaccessible random time $T^i$ fulfills the property 
	$$ 
	\one_{[[T^i]]}(\omega, S(\omega))\,\one_{\{S< \infty\}}= 0,
	$$ 
	for every predictable random time $S$, 	
	see  Definition 2.20, Chapter I,  in \cite{JacodBook}.	\end{itemize} 
\end{remark} 

We recall an important notion of measure associated with  $\mu$, given in formula (3.10) in \cite{jacod_book}. 
\begin{definition} 
	\label{D_DoleansMeas} 
	Let $(\tilde \Omega_n)$ be a  partition  of $\tilde \Omega$ constituted by elements of $\tilde{\mathcal O}$, such that $\one_{\tilde \Omega_n} \star \mu \in \mathcal A$. 
	$M^{\P}_{\mu}$ denotes the $\sigma$-finite measure on  $(\tilde{\Omega}, 
	\mathcal{{\tilde F}^{\P}})$, such that 
	for every $W: \tilde \Omega \rightarrow \R $ positive, bounded,   $\mathcal{\tilde{F}}^{\P}$-measurable function, 
	\begin{equation}\label{DoleansMeas} 
	M^{\P}_{\mu}(W\,\one_{\tilde{\Omega}_n})= \E\big[W\,\one_{\tilde{\Omega}_n} \star \mu_T\big]. 
	\end{equation} 	 
\end{definition} 

Let us now  set $\hat{\nu}_t(de):= \nu(\{t\}, de)$ for all $t \in [0,\,T]$. 
For 
any 
$W \in \mathcal{\tilde{O}}$, 
we define 
\begin{align*} 
	\hat{W}_t = \int_{\R} W_t(x)\,\hat{\nu}_t(de), 
	\quad \tilde{W}_t = \int_{\R} W_t(x)\,\mu(\{t\}, de)- \hat{W}_t,\quad t \geq 0,\nonumber
\end{align*}
with the convention that
	$\tilde{W}_t= + \infty$ if $\hat{W}_t$ is not defined. 
For every $q \in [1,\,\infty[$, we  introduce the linear spaces 
\begin{align*}
	\mathcal{G}^q(\mu)=\Big\{W \in \mathcal{\tilde{P}} 
	:\,\, 
	\Big[\sum_{s\leq \cdot}|\tilde{W}_s|^2\Big]^{q/2}\in \mathcal{A}^+\Big\},\quad 
	\mathcal{G}^q_{\textup{loc}}(\mu)=\Big\{W \in  \mathcal{\tilde{P}} 
	: \,\, 
	\Big[\sum_{s\leq \cdot}|\tilde{W}_s|^2\Big]^{q/2}\in \mathcal{A}_{\textup{loc}}^+\Big\}.\nonumber
\end{align*} 
Given $W \in \mathcal{\tilde{P}}$, we define 
the  increasing (possibly infinite) predictable process
\begin{align} 
	C(W)_t &= |W - \hat{W}|^2\star \nu_t + \sum_{s \leq t}(1-\hat{\nu}_s(\R))\,|\hat{W}_s|^2,\label{C(W)}
\end{align} 
provided the right-hand side is well-defined.
By  
Proposition 3.71 in \cite{jacod_book}, we have
$
	\mathcal{G}^2(\mu)=\{ W \in \mathcal {\tilde{P}}:\,\, 
||W||_{\mathcal{G}^2(\mu)}< \infty\}
$, where 
\begin{align} 
||W||^2_{\mathcal{G}^2(\mu)}&:=\sper{C(W)_T}.
	\label{G2mu_norm} 
\end{align} 
We also introduce the  space 
\begin{equation}\label{L2mu} 
	\mathcal{L}^2(\mu):=\left\{W \in \tilde{\mathcal{P}} 
	:\,\, ||W||_{\mathcal{L}^2(\mu)}< \infty\}, \,\, \textup{with}\,\,\, ||W||_{\mathcal{L}^2(\mu)}:=\E\Big[\int_{]0,T]\times \R}  |W_s(e)|^2 \,\nu(ds\,de)\Big]\right\}.
\end{equation} 

\begin{lemma}\label{L_G2mu_L2mu_inclusion} 
	Let  $p=1,2$. If $|W|^p \star \mu \in \mathcal{A^+_{\rm loc}}$ then $W \in \mathcal{G}^p_{\rm loc}(\mu)$. 
\end{lemma}
\proof
The case $p=1$ is stated in Proposition 1.28, Chapter II,  in \cite{JacodBook}.
In order to prove  the case $p=2$,
 it is enough to show that, if $W \in \mathcal{L}^2(\mu)$, then   
$ 
		||W||^2_{\mathcal{G}^2(\mu)} \leq ||W||^2_{\mathcal{L}^2(\mu)}$.
Indeed, if this holds, then the conclusion would follow by usual localization arguments.

Let $W \in \mathcal{\tilde{P}}$. 
Recalling \eqref{C(W)} and \eqref{G2mu_norm}, we evaluate the expectation of the right-hand side in \eqref{C(W)}.
For every $t \geq 0$, since $\hat \nu_t(\R) \leq 1$, we have 
	\begin{equation} 
	\sum_{s \in ]0,\,t]}|\hat{W}_s|^2(1-\hat \nu_s(\R)) \leq \sum_{s \leq t}|\hat{W}_s|^2 
	\leq \sum_{s \leq t}\hat \nu_s(\R)\int_{\R}|W_s(e)|^2\,\hat{\nu}_s(de)\leq |W|^2 \star \nu_t.\label{ineq_sum_L2mu}
	\end{equation} 
Moreover, taking into account that 
$|\hat W|^2 \star \nu_t= \sum_{s \leq t}|\hat{W}_s|^2\,\hat{\nu}_s(\R)$, for every $t \geq 0$,
the process  $C(W)$ defined in \eqref{C(W)} 
can be decomposed as 
\begin{align}\label{C(W)_dec} 
C(W)_t&= |W|^2 \star \nu_t - 2\sum_{s \leq t} |\hat{W}_s|^2+ \sum_{s \leq t}|\hat{W}_s|^2\,\hat{\nu}_s(\R) + \sum_{s \leq t}|\hat{W}_s|^2\,(1-\hat{\nu}_s(\R))\nonumber\\ 
&=	|W|^2 \star \nu_t - \sum_{s \leq t} |\hat{W}_s|^2. 
\end{align} 
In particular, since $W \in \mathcal{L}^2(\mu)$, we have 
\[ 
||W||^2_{\mathcal{G}^2(\mu)} =  \E\Big[\int_{]0,T]\times \R}|W_s(e)|^2\,\nu(ds\,de) - \sum_{s \in ]0,T]} |\hat{W}_s|^2\Big]\leq||W||^2_{\mathcal{L}^2(\mu)}. 
\] 
This concludes the proof.
\endproof

\begin{proposition}\label{Sec:WD_R_varphi_integr} 
	Let $\varphi:\Omega \times [0,T] \times \R \rightarrow \R$ be a $\mathcal{\tilde P}$-measurable function and  $A$ a $\mathcal {\tilde P}$-measurable subset of $\Omega \times [0,\,T] \times \R$,   such that 
	\begin{align} 
	&|\varphi|\,\one_{A} \star \mu^X \in \mathcal A^+_{\rm loc},\label{L1int}\\ 
	&|\varphi|^2\,\one_{A^c} \star \mu^X \in \mathcal A^+_{\rm loc}\label{G2int}. 
	\end{align} 
	Then the process $\varphi$ belongs to $\mathcal{G}^1_{\rm loc}(\mu^X)$. 
\end{proposition}
\proof	 
Formula	\eqref{L1int}, together with 
Lemma \ref{L_G2mu_L2mu_inclusion}  with $p=1$, 
 gives that $\varphi\,\one_{A}$ belongs to $\mathcal G^1_{\rm loc}(\mu^X)$. On the other hand,  
%
by Lemma \ref{L_G2mu_L2mu_inclusion} with $p=2$, formula	\eqref{G2int} implies that 	$\varphi\,\one_{A^c}$ belongs to $\mathcal G^2_{\rm loc}(\mu^X) \subset \mathcal G^1_{\rm loc}(\mu^X)$. 
\endproof


\begin{remark}\label{R_FARE} 
Since
$\hat{W}= \hat{W}\,\one_{J}$,
$\hat{\nu}(\R)\,\one_{K}=\one_{K}$,
$1-\hat{\nu}(\R)>0$ on 
$J \setminus K$, 
the quantity $C(W)$ in \eqref{C(W)} can be rewritten as 
\begin{align}
C(W) = |W - \hat{W}\,\one_{J}|^2\star \nu + \sum_{s \leq \cdot}(1-\hat{\nu}_s(\R))\,|\hat{W}_s|^2\,\one_{J \setminus K}(s).\label{C(W)BIS}
\end{align}
\end{remark}

\begin{proposition}	\label{P_J=K_totally_inaccessible} 
	If $D$ is the disjoint union of $K$ and  $\cup_n [[T^i_n]]$, where $(T^i_n)_n$ are totally inaccessible finite times,  then  $J=K$ up to an evanescent set. 
\end{proposition} 
\proof 
We start by noticing some basic facts. 
By Remark \ref{R_pred_supp}, we have 
$\one_J= 
{}^p{}(\one_D)$,
and
${}^p{} 
(\one_K ) 
= 
\one_K$; 
on the other hand, by Theorem 5.2  in \cite{chineseBook}, 
for any  totally inaccessible time $T^i$ we have 
that ${}^{p}\left(\one_{[[T^i]]}\,\one_{\{T^i < \infty\}}\right)$ $=0$,
and therefore 
the predictable projection of $\one_{[[T^i_n]]}$ is zero  
since $T^i_n$ is a 
totally inaccessible  finite time. 
Consequently, by additivity of the predictable projections,  we obtain 
${}^p{}(\one_D)={}^p{}(\one_K) + \sum_n{}^p{}(\one_{[[T^i_n]]})={}^p{}(\one_K)$.
Collecting previous identities,
we conclude that
$ 
\one_J= {}^p{}(\one_D)= \one_K, 
$
therefore 
$ 
J=K. 
$ 
\endproof

 \begin{proposition}\label{P_forBSDEs2} 
If 
$C(W)_T=0$ a.s.,
then 
\begin{equation}\label{W-hatW_K} 
||W-\hat{W}\,\one_{K}||_{\mathcal{L}^2(\mu)}=0,
 \end{equation} 
or, equivalently,  there exists a predictable process $(l_s)$ such that
\begin{equation}\label{2.12}
W_s(e)= l_s\,\one_{K}(s),\quad d\P\,\nu(ds\,de)\textup{-a.e.}
\end{equation}
In particular, 
\begin{equation}\label{nu_c_new}
W_s(e)=0, \quad d\P\,\nu^c(ds\,de)\textup{-a.e.}, 
\end{equation}
and  there is a predictable process  $(l_s)$  such that 
\begin{equation}\label{nu_d_new}
W_s(e)= l_s\,\one_{K}(s), \quad d\P\,\nu^d(ds\,de)\textup{-a.e.}	
\end{equation}
\end{proposition} 
\proof 
By \eqref{C(W)BIS} we have
\begin{equation*} 
\begin{cases} 
	|W - \hat{W}\,\one_{J}|^2\star \nu =0, \\ 
	\sum_{s \leq \cdot}(1-\hat{\nu}_s(\R))\,|\hat{W}_s|^2\,\one_{J \setminus K}(s)=0. 
\end{cases} 
\end{equation*} 
Since $1-\hat{\nu}(E)>0$ on $J \setminus K$,
previous identities  give 
\begin{equation*} 
\begin{cases} 
|W - \hat{W}\,\one_{J}|^2\star \nu =0, \\ 
\hat{W}\,\one_{J \setminus K}=0, 
\end{cases} 
\end{equation*} 
and this  gives \eqref{W-hatW_K}.
We show now the equivalence property. Obviously \eqref{W-hatW_K}  implies \eqref{2.12} setting $l = \hat W$. Conversely, for fixed $s$ such that $(\omega, s) \in K$, integrating \eqref{2.12} against $\one_{\{s\} \times \R}$, we obtain $\hat{W_s} \one_K(s)= l_s \one_K(s)$, which implies \eqref{W-hatW_K}.
Finally, \eqref{nu_c_new} and \eqref{nu_d_new} follow observing that  $|W - l\,\one_{K}|^2\star \nu = |W|^2\star\nu^c+|W-l\,\one_{K}|^2\star\nu^d$.
\endproof 

 \subsection{A class of stochastic processes $X$ related in a specific way to an integer-valued random measure $\mu$}\label{Sec_2.2}
We will formulate  two assumptions related to an integer-valued random measure $\mu$ on $[0,\,T] \times \R$ and some 
 c\`adl\`ag process $X$. 
 We recall that a sequence of random times  $(T_n)_n$ exhausts the jumps of a process $Y$ if  
		$[[T_n]] \cap [[T_m]] = \emptyset$, $n \neq m$, and $\{\Delta Y \neq 0 \} = \cup_n [[T_n]]$, see Definition 1.30, Chapter I, in \cite{JacodBook}.
A process $Y$ is quasi-left-continuous if and only if there is a sequence of totally inaccessible times  $(T^i_n)$ that exhausts the jumps of $Y$, see Proposition 2.26, Chapter I,  
in \cite{JacodBook}.
		

\begin{hypothesis}\label{H_X_mu} 
	$X$ is an adapted c\`adl\`ag process  with decomposition $X=X^i + X^p$, where
\begin{enumerate}
	\item $Y:=X^i$ is a c\`adl\`ag quasi-left-continuous adapted process satisfying	$\{\Delta Y \neq 0\} \subset D$. 
	Moreover, there exists a $\tilde{\mathcal{P}}$-measurable map  $\tilde \gamma: \Omega \times ]0,\,T] \times \R\rightarrow \R$ such that 
	\begin{equation}\label{B35ii} 
	\Delta Y_t(\omega)\,\one_{]0,\,T]}(t)=  \tilde \gamma(\omega,t,\cdot)\quad dM^\P_{\mu}\textup{-a.e.} 
	\end{equation} 
  \item $X^p$ is a c\`adl\`ag predictable process satisfying 	$\{\Delta X^p \neq 0\} \subset J$. 
\end{enumerate}
\end{hypothesis}
\begin{remark}\label{R_HP_Jacod} 
	Theorem 3.89 in \cite{jacod_book} states an 
	It\^o formula which transforms a 
	special semimartingale $X$ into a special semimartingale $F(X_t)$ through 
	a  $C^{2}$ function $F: \R \rightarrow \R$. 
	There, the process $Y=X$ is supposed to fulfill 	Hypothesis \ref{H_X_mu}-1.
\end{remark} 
\begin{remark} 
	If $\mu$ is the jump measure of a c\`adl\`ag process $X$, then  
Hypothesis 
	\ref{H_X_mu}-1. holds for $Y= X$,	
	with $\tilde \gamma(t,\omega, x)=x$. 
\end{remark} 
 
The proof of the following result is reported in Section \ref{Sec_Proof_P_identiy_measures}. 
\begin{proposition}\label{P_identiy_measures} 
	Let  $X$ be a process verifying Hypothesis \ref{H_X_mu}.
	Then, there exists a null set $\mathcal N$ such that, 
	for every Borel function $\varphi: [0,\,T] \times  \R \rightarrow \R_+$ satisfying  $\varphi(s,0)=0$, $s \in [0,\,T]$, we have, for every $\omega \notin \mathcal N$, 
	\begin{equation}\label{E_claim} 
	\int_{]0,\,T]\times \R}\varphi(s,x)\,\mu^X(\omega,ds\,dx)= \int_{]0,\,T]\times \R}\varphi(s,\tilde \gamma(\omega,s,e))\,\mu(\omega,ds\,de) + V^{\varphi}(\omega), 
	\end{equation} 
	with  $V^{\varphi}(\omega)=\sum_{0<s \leq T} \varphi(s,\Delta X^p_s(\omega))$. 
	In particular, 
	\begin{equation}\label{E_claim2} 
	\int_{]0,\,T]\times \R}\varphi(s,x)\,\mu^X(\omega,ds\,dx)\geq \int_{]0,\,T]\times \R}\varphi(s,\tilde \gamma(\omega,s,e))\,\mu(\omega,ds\,de) \quad \textup{for every}  \,\,\omega \notin \mathcal N. 
	\end{equation} 
	Identity \eqref{E_claim} still holds true when $\varphi:[0,\,T] \times \R \rightarrow \R$ and the left-hand side is finite. 
\end{proposition} 
\begin{remark}\label{R_phi_real_val}
	The result in Proposition \ref{P_identiy_measures} still holds true if $\varphi$ is a real-valued random  function  on $\Omega \times [0,\,T]\times \R$. 
\end{remark}

We will make use in the sequel of the following assumption on  $\mu$. 
\begin{hypothesis}\label{H_nu} 
	\hspace{2em} 
	\begin{itemize} 
		\item[(i)]$D= K \cup (\cup_n [[T^i_n]])$ up to an evanescent set, where $(T^i_n)_n$ are totally inaccessible times such that 
		$[[T^i_n]] \cap [[T^i_m]] = \emptyset$, $n \neq m$; 
		\item[(ii)]for every predictable time $S$ such that $[[S]] \subset K$, 
		$\nu(\{S\},de)=\mu(\{S\},de)$ a.s. 
	\end{itemize} 
\end{hypothesis} 
\begin{remark}\label{R_Hp_3.1} 
	Hypothesis \ref{H_nu}-(i)  implies that 
	$J=K$, up to an evanescent set, see Proposition \ref{P_J=K_totally_inaccessible}. 
\end{remark} 
\begin{remark}\label{R_H_mu_BIS} 
\begin{enumerate}
\item	If $\hat \nu =0$ then 	 
		 $J = K = \emptyset$. Taking into account  Remark \ref{R_H_mu}-(ii),  $D=   \cup_n [[T^i_n]]$, $[[T^i_n]] \cap [[T^i_m]] = \emptyset$, $n \neq m$, where $(T^i_n)_n$ are totally inaccessible times, and 
		Hypothesis 
		\ref{H_nu} trivially holds. 
\item Notice that, if 
  $\nu$ is given in the form \eqref{nu_dis} 
and the process $A$ appearing in  \eqref{nu_dis} is continuous, then $\hat \nu =0$.
 In that case, Theorem (4.47) of \cite{jacod_book} states
that $\mu$ is quasi-left-continuous.
\end{enumerate}
\end{remark} 
 
We have the following important technical result, for the proof see Section \ref{Sec_Proof_L_ident_mu_muX}.
\begin{proposition}\label{L_ident_mu_muX} 
	Let  $\mu$ satisfy Hypothesis \ref{H_nu}. 
	Assume that  $X$ is a process verifying  Hypothesis \ref{H_X_mu}. 
	Let  $\varphi: \Omega \times [0,\,T] \times \R \rightarrow \R_+$ 
	such that $\varphi(\omega, s,0)=0$ for every $s \in [0,\,T]$, up to indistinguishability, and assume that there exists  a $\mathcal {\tilde P}$-measurable subset $A$ of $\Omega \times [0,\,T] \times \R$ satisfying 
	\begin{equation}\label{cond_int_A} 
	|\varphi|\,\one_{A} \star \mu^X \in \mathcal A^+_{\rm loc},\quad |\varphi|^2\,\one_{A^c} \star \mu^X \in \mathcal A^+_{\rm loc}. 
	\end{equation} 
	Then 
	\begin{equation}\label{Id_intstoch} 
	\int_{]0,\,\cdot]\times \R}\varphi(s,x)\,\,(\mu^X-\nu^X)(ds\,dx) 
	= \int_{]0,\,\cdot]\times \R}\varphi(s,\tilde \gamma(s,e))\,(\mu-\nu)(ds\,de).
	\end{equation} 
\end{proposition}

\begin{remark} 
	Under  condition \eqref{cond_int_A}, Proposition \ref{Sec:WD_R_varphi_integr} and  inequality \eqref{E_claim2} in Proposition \ref{P_identiy_measures} imply  that 
	$(s,x) \mapsto \varphi(s,x) \in \mathcal G^1_{\rm loc}(\mu^X)$ and 
	$(s,e) \mapsto \varphi(s,\tilde \gamma(s,e)) \in \mathcal G^1_{\rm loc}(\mu)$. 
	In particular the two stochastic integrals in \eqref{Id_intstoch} are well-defined. 
\end{remark}

We end the section focusing on the case when $X$ is of jump-diffusion type. The following result is proved in Section \ref{Sec_Proof_Ex_guida}.
\begin{lemma}\label{Ex_guida} 
	Let $\mu$ 
	satisfy Hypothesis 
	\ref{H_nu}. 
	Let 
	$N$ be a continuous martingale, and $B$ an increasing predictable c\`adl\`ag  process, with $B_0=0$, such that 
	$\{\Delta B \neq 0\} 
	\subset J$. 
	Let $X$ be a process  which is solution of equation 
	\begin{equation}\label{X_SDE} 
	X_t = X_0 +\int_0^t b(s,X_{s-})\,d B_s +\int_0^t \sigma(s,X_s)\,dN_s + \int_{]0,\,t]\times \R} \gamma(s,X_{s-},e)\,(\mu-\nu)(ds\,de), 
	\end{equation} 
	for some given Borel functions $b,\sigma: [0,\,T]\times \R \rightarrow \R$, and $\gamma:[0,\,T]\times \R \times \R \rightarrow \R$ such that 
	\begin{align} 
	&\int_0^t |b(s,X_{s-})|\, d B_s < \infty \,\,\textup{a.s.},\label{Ver_1}\\ 
	&\int_0^t |\sigma(s,X_s)|^2\,d [N,N]_s < \infty\,\, \textup{a.s.},\label{Ver_2}\\ 
	&(\omega,s,e)\mapsto \gamma(s,X_{s-}(\omega),e)\in \mathcal G^1_{\rm loc}(\mu).\label{Ver_3} 
	\end{align}	 
Then $X$ satisfies  Hypothesis \ref{H_X_mu} with decomposition $X=X^i+X^p$, where
	\begin{align} 
	X^i_t &= \int_{]0,\,t]\times \R} \gamma(s,X_{s-},e)\,(\mu-\nu)(ds\,de),\label{Xi}\\ 
	X^p_t &= X_0 +\int_0^t b(s,X_{s-})\,d B_s +\int_0^t \sigma(s,X_s)\,dN_s,\label{Xp} \\
	\tilde \gamma(\omega, s,e) &= \gamma(s,X_{s-}(\omega),e)\,(1-\one_{K}(\omega, s)).\nonumber
	\end{align} 

\end{lemma}

 
\section{BSDEs driven by an integer-valued random measure} \label{Sec_application_BSDE} 
Let $\mu$ be an integer-valued random measure defined on $[0,T]\times \R$. 
Let $M$ be a continuous process with $M_0=0$. 
Let $(\mathcal F_t)$ be the canonical filtration associated to $\mu$ and $M$, and suppose that $M$ is an $(\mathcal F_t)$-local martingale. 
$\nu$  will denote a "good" version of the dual predictable projection of $\mu$   in the sense of 
Proposition 1.17, Chapter II, in \cite{JacodBook}. In particular, $\nu(\omega, \{t\}\times \R) \leq 1$ identically. 
Let $\lambda$ be  a predictable random measure on $[0,\,T]\times \R$, and  
$\zeta$  a non-decreasing adapted continuous process.

 We focus now on the general BSDE \eqref{GeneralBSDE}
of the Introduction, whose coefficients are the following: 
$\xi$ is an $\mathcal F_T$-measurable  square integrable random variable, 
$\tilde f: \Omega \times [0,\,T] \times \R^3 \rightarrow \R$ (resp. $\tilde g: \Omega \times [0,\,T]\times \R^2\rightarrow \R$)
 is a measurable function, whose domain is equipped with the $\sigma$-field 
 $\mathcal F\otimes \mathcal B([0,\,T] \times \R^3)$ (resp. $\mathcal F\otimes \mathcal B([0,\,T] \times \R^2)$).
 
A solution  of BSDE \eqref{GeneralBSDE} is a triple of processes $(Y,Z,U)$ such that  the first two  integrals in \eqref{GeneralBSDE} exist and are finite in the Lebesgue sense, $Y$ is adapted and c\`adl\`ag, $Z$ is  progressively measurable with  $Z \in L^2([0,\,T], d \langle M\rangle_t)$ a.s., and  
   $U \in \mathcal{G}^2_{\rm loc}(\mu)$.

\begin{remark}\label{R_uniq_G2}
Uniqueness in $\mathcal G^2_{\rm loc}(\mu)$ 
means the following: 
if $(Y,Z,U)$, $(Y',Z',U')$ are solutions of the BSDE 
\eqref{GeneralBSDE}, then $Y=Y'$ $\lambda(ds,\R)$ and $d \zeta_s$-a.e., for almost all $\omega$, $Z=Z'$ $d \P\,d \langle M\rangle_t$ a.e., and
there is a predictable process $(l_t)$ such that 
$U_t(e)-U'_t(e) = l_t\,\one_K(t)$, $d\P\,\nu(dt\,de)$-a.e.	
The latter fact is a direct consequence of Proposition \ref{P_forBSDEs2}. 
Moreover, provided that $\lambda <\!< \nu$, given a solution  $(Y,Z,U_0)$ of BSDE \eqref{GeneralBSDE}, 
the class of all solutions will be given 
by the  triples $(Y,Z,U)$, where 
$U= l\,\one_K + U_0$ for some predictable process  $(l_t)$.	 
In particular, if  $K = \emptyset$, 
then the third component  of the  BSDE solution is 
uniquely characterized in  
  $\mathcal{L}^2(\mu)$. 
\end{remark}
\begin{remark}\label{Rem_well-posedness:Xia} 
	A general BSDE of type \eqref{GeneralBSDE} is considered for instance in \cite{xia}, 
	with the following restrictions on the random measures $\lambda$ and $\nu$: 
	\begin{align} 
	&\lambda([0,\,T]\times \R)\,\,\textup{ is a bounded random variable,  }\,\lambda([0,\,t]\times \R)\,\,\textup{ is continuous w.r.t. }\,\,t,\nonumber\\ 
	&\nu([0,\,t]\times \R)\,\,\textup{ is continuous w.r.t. }\,t.\label{nu_cont} 
	\end{align} 
	In Theorem 3.2 in \cite{xia}, the author proves  that  under suitable assumptions on the coefficients $(\xi,\tilde f, \tilde g)$ 
	there exists  a unique  triplet of processes $(Y,Z,U) \in \mathcal{L}^2(\zeta\lambda) \times \mathcal{L}^2(M) \times \mathcal L^2(\mu)$,  with 
	$\E\big[\sup_{t \in [0,\,T]}Y_t^2 \big] < \infty$, satisfying BSDE \eqref{GeneralBSDE}, where 
	\begin{align*} 
	\mathcal{L}^2(\zeta\lambda):&=\Big\{(Y_t)_{t \in [0,\,T]}\,\, \textup{optional  :}\,\,\E\Big[\int_0^T Y^2_s \,d \zeta_s\Big] + \E\Big[\int_0^T Y^2_s \,\lambda(ds,\R)\Big] < \infty\Big\},\nonumber\\ 
	\mathcal{L}^2(M):&=\Big\{(Z_t)_{t \in [0,\,T]}\,\,\textup{predictable  :}\,\,\E\Big[\int_0^T Z^2_s \,d \langle M \rangle_s\Big]  < \infty\Big\},\nonumber 
	\end{align*}	 
	and $\mathcal{L}^2(\mu)$ is the space introduced in \eqref{L2mu}. 
\end{remark} 

When  $\zeta$ and $M$ vanish, BSDE  \eqref{GeneralBSDE} turns out to be driven only by a purely discontinuous martingale, and  becomes 
\begin{align}\label{GeneralBSDE_disc} 
Y_t &= \xi + \int_{]t,\,T]\times \R} \tilde f(s,\,e,\,Y_{s-},\,U_{s}(e))\,  \lambda(ds\,de) - \int_{]t,\,T]\times \R} U_s(e)\,(\mu-\nu)(ds\,de). 
\end{align} 
\begin{remark}\label{R_orth_BSDE}
	The process $Y$ solution to \eqref{GeneralBSDE_disc}  is an $(\mathcal F_t)$-orthogonal process. In fact, for every continuous $(\mathcal F_t)$-local  martingale $N$ we have 
\begin{equation}\label{YNcov} 
[Y,N]_t = \left[\int_{]0,\,\cdot]\times \R} \tilde f(s,\,e,\,Y_{s-},\,U_{s}(e))\,  \lambda(ds,de), N\right]_t - 
\left[\int_{]0,\,\cdot]\times \R} U_s(e)\,(\mu-\nu)(ds\,de), N\right]_t. 
\end{equation} 
Since $\int_{]0,\,\cdot]\times \R} \tilde f(s,\,e,\,Y_{s-},\,U_{s}(e))\,  \lambda(ds,de)$ is a  bounded variation process, the first bracket in the right-hand side of \eqref{YNcov} is zero  by Proposition 3.14 in \cite{BandiniRusso1}. 
On the other hand, the second term in the right-hand side of \eqref{YNcov} is zero because  $\int_{]0,\,\cdot]\times \R} U_s(e)\,(\mu-\nu)(ds\,de)$ is a purely discontinuous martingale. 
\end{remark}


 
\subsection{Identification of the BSDEs solution}\label{Section_BSDEs_applications}


The fundamental tool of this section is the notion of covariation of two processes $X$ and $Y$ (that are not necessarily semimartingales), denoted $[X,Y]$,  see Definition 3.4 in \cite{BandiniRusso1}. A process $X$ is said to be finite quadratic variation process if $[X,X]$ exists.
Any $({\mathcal F}_t)$-adapted process  $X$  is said to be $({\mathcal F}_t)$-orthogonal if 
	$[X, N] = 0$  for every $N$ continuous local $({\mathcal F}_t)$-martingale. 

	
By Proposition 5.8 in \cite{BandiniRusso1} we have the following result.
	\begin{proposition}\label{Sec:WD_P_unique_decomp}
	Any $(\mathcal F_t)$-special weak Dirichlet process $X$ admits a 
 decomposition of the type
 \begin{equation}\label{Mc+Md+A} 
	X=X^c+M^d+A, 
	\end{equation} 
	where $X^c$ is a continuous local martingale, $M^d$ is a purely discontinuous local martingale, and $A$ is an $({\mathcal F}_t)$-predictable and orthogonal process, with $A_0=0$.
	\eqref{Mc+Md+A} is called the canonical decomposition of $X$.
\end{proposition}


 



 
The following condition on  $X$ will play a fundamental role in the sequel: 
\begin{equation}\label{Sec:WD_CNS} 
\int_{]0,\cdot]\times \R} |x|\,\one_{\{|x| >1\}}\,\mu^X(ds\,dx) \in \mathcal{A}^{+}_{\rm loc}. 
\end{equation} 
%
%
%
Moreover, we will be interested in functions  $v:[0,T] \times \R \rightarrow \R$  fulfilling the  integrability property \begin{equation}\label{Sec:WD_E_C} 
\int_{]0,\cdot]\times \R} |v(s,X_{s-} +x )-v(s,X_{s-})-x\,\partial_x v(s,X_{s-})|\,\one_{\{|x| >1\}}\,\mu^X(ds\,dx) \in \mathcal{A}^{+}_{\rm loc}. 
\end{equation}



\begin{remark}
By Proposition 2.7 in \cite{BandiniRusso1}, if $X$ is a c\`adl\`ag  process  satisfying 
	condition \eqref{Sec:WD_CNS}, and     $v:[0,T] \times \R \rightarrow \R$  is  a function of class $C^{0,1}$, then	 
	\begin{align} 
		|v(s,X_{s-} + x)-v(s,X_{s-})|^2\,\one_{\{\vert x \vert \leq 1\}}\star \mu^X \in \mathcal{A}^+_{\rm loc}.\label{A2_Aloc_NEW} 
		\end{align}	 
Moreover,  Lemma 5.24 in \cite{BandiniRusso1} states  that
if $X$ is a c\`adl\`ag  process  satisfying 
	condition \eqref{Sec:WD_CNS}, and     $v:[0,T] \times \R \rightarrow \R$  is  a function of class $C^{0,1}$ fulfilling \eqref{Sec:WD_E_C}, then	 
\begin{align} 
		|v(s,X_{s-} + x)-v(s,X_{s-})|\,\one_{\{\vert x \vert > 1\}}\star \mu^X \in \mathcal{A}^+_{\rm loc}.\label{F2_Aloc_NEW} 
		\end{align}	 
\end{remark}




We have the following result.
\begin{proposition}\label{R_gammatilde_id} 
	Let  $X$ be a process such that $(X,\mu)$ verifies Hypothesis \ref{H_X_mu}. 
	Let in addition $v:[0,T] \times \R \rightarrow \R$ be a function of class $C^{0,1}$. 
		If  $X$ and $v$ satisfy  conditions \eqref{Sec:WD_CNS} and   \eqref{Sec:WD_E_C}, and  moreover  $\sum_{s \leq T}|\Delta X_{s}|^2 < \infty $ a.s., 
		then 
		$$ 
		(s,e) \mapsto v(s,X_{s-} +\tilde{\gamma}(s,e) )-v(s,X_{s-}) \in \mathcal G^1_{\rm loc}(\mu). 
		$$ 
\end{proposition} 
\proof
\emph{Step 1.}  We first notice that, if	 
		$\sum_{s \leq T}|\Delta X_{s}|^2 < \infty $ a.s., 
		then 
		\begin{align*}
		|v(s,X_{s-} + \tilde \gamma(s,e))-v(s,X_{s-})|^2\,\one_{\{\vert \tilde \gamma(s,e) \vert \leq 1\}}\star \mu \in \mathcal{A}^+_{\rm loc}.
		\end{align*}	 
Indeed, this follows from \eqref{A2_Aloc_NEW} and inequality \eqref{E_claim2} in  Proposition \ref{P_identiy_measures}, 
 with 
$\varphi(\omega, s,x)=|v(s,X_{s-}(\omega) + x)-v(s,X_{s-}(\omega))|^2\,\one_{\{\vert  x \vert \leq 1\}}$, taking into account  Remark \ref{R_phi_real_val}.

\noindent \emph{Step 2.} We observe that, if $X$ and $v$ satisfy  conditions \eqref{Sec:WD_CNS} and   \eqref{Sec:WD_E_C}, 
		then 
		\begin{align*} 
		|v(s,X_{s-} +\tilde{\gamma}(s,e) )-v(s,X_{s-})|\,\one_{\{|\tilde{\gamma}(s,e)| >1\}}\star \mu \in \mathcal{A}^{+}_{\rm loc}.
		\end{align*} 
This is a consequence  \eqref{F2_Aloc_NEW}, 
\eqref{Sec:WD_CNS} 
and \eqref{Sec:WD_E_C} 
together with 
inequality \eqref{E_claim2} in Proposition \ref{P_identiy_measures}, with  
$\varphi(\omega,s,x)=  |v(s,X_{s-}(\omega) +x )-v(s,X_{s-}(\omega))|$ $\one_{\{|x| >1\}}$, 
taking into account  Remark \ref{R_phi_real_val}.
 
\noindent \emph{Step 3.} 
The conclusion of Proposition \ref{R_gammatilde_id} is a direct consequence of Steps 1 and 2, 
together with Proposition  \ref{Sec:WD_R_varphi_integr}, with $\varphi(\omega, s,e)=v(s,X_{s-}(\omega) + \tilde \gamma(\omega, s,e))-v(s,X_{s-}(\omega))$ and $A = \{(\omega,s,e): |\tilde\gamma(\omega, s,e)| >1 \}$. 
\endproof

Let us now consider the following  assumption on a couple $(X,Y)$ of adapted processes.
\begin{hypothesis}\label{H_chain_rule_C01} 
	$X$ is a  special weak Dirichlet process      of finite quadratic variation, satisfying condition \eqref{Sec:WD_CNS}. 
	$Y_t = v(t,\,X_t)$ for some (deterministic) function $v: [0,T] \times \R \rightarrow \R$ of class $C^{0,1}$ such that $v$ and $X$ verify   condition \eqref{Sec:WD_E_C}. 
\end{hypothesis} 
Our first main result is the following. 
\begin{theorem}\label{P_ident} 
	Let $\mu$ 
	satisfy Hypothesis \ref{H_nu}, and assume that 
	 $X$ is a process such that $(X,\mu)$ verifies  Hypothesis \ref{H_X_mu}. 
	Let  $(Y,Z,U)$ be a solution to the BSDE \eqref{GeneralBSDE} such that the pair $(X, Y)$ satisfies 
	Hypothesis \ref{H_chain_rule_C01} with corresponding  function $v$. 
Let $X^c$ denote the continuous local martingale of $X$ given in the canonical decomposition \eqref{Mc+Md+A}. 
	 
	Then,  the pair $(Z,U)$ 
	fulfills 
	\begin{equation}\label{id_2BIS} 
	Z_t = \partial_x v(t,X_t) \,\frac{d\langle X^c,M\rangle_t}{d\langle M \rangle_t}\quad 
	d\P \,d\langle M \rangle_t \,\textup{-a.e.,} 
	\end{equation} 
	\begin{equation}\label{id_3} 
	\int_{]0,\,t]\times \R} H_s(e)\,(\mu-\nu)(ds\,de)=0,\quad  \forall \,\, t \in ]0,\,T],\,\, \textup{a.s.,} 
	\end{equation} 
	with	 
	\begin{equation}\label{K_def} 
	H_s(e):=U_s(e)-(v(s,X_{s-}+ \tilde\gamma(s,e))-v(s,X_{s-})). 
	\end{equation} 
	If, in addition, $H \in \mathcal G^2_{\rm loc}(\mu)$, 
	then there exists a predictable process $(l_s)$ such that 
	$$
	H_s(e)= l_s\,\one_{K}(s),\quad d\P\,\nu(ds\,de)\textup{-a.e.}
	$$ 
	In particular, 
\begin{equation}\label{ID_cont_1}
H_s(e)= 0, \quad d\P\,\nu^c(ds\,de)\textup{-a.e.}
\end{equation}
and	
\begin{equation}\label{ID_disc_1}
H_s(e)= l_s\,\one_{K}(s), \quad d\P\,\nu^d(ds\,de)\textup{-a.e.}
\end{equation}
\end{theorem} 
\begin{remark} 
Notice that $H$ in  \eqref{K_def} belongs to  $\mathcal G^1_{\rm loc}(\mu)$, so that the integral in \eqref{id_3} is well-defined. 
Indeed, 
by Hypothesis \ref{H_chain_rule_C01},  
  $X$ and $v$ in the statement of Theorem \ref{P_ident} satisfy \eqref{Sec:WD_CNS} and \eqref{Sec:WD_E_C}. 
	By Proposition \ref{R_gammatilde_id} 
	it follows that  $(s,e) \mapsto (v(s, X_{s-}+ \tilde \gamma(s,e))-v(s, X_{s-})) \in \mathcal G^1_{\rm loc}(\mu)$. 
	Since  $U \in \mathcal G^2_{\rm loc}(\mu) \subset  \mathcal G^1_{\rm loc}(\mu)$, the conclusion follows.
\end{remark}

The proof of Theorem \ref{P_ident} is based essentially on the following stability result for 
c\`adl\`ag processes,  which was   the object  of Theorem 5.26  in \cite{BandiniRusso1}. 
\begin{theorem}\label{T_C1_dec_specialweak_Dir} 
		Let $X$ be an $(\mathcal F_t)$-special weak Dirichlet process  of finite quadratic variation  with its canonical decomposition  $X=X^c+M^d+A$, satisfying 	condition \eqref{Sec:WD_CNS}. 
		Then,  
		for every  $v:[0,T] \times \R \rightarrow \R$ of class  $C^{0,1}$ verifying \eqref{Sec:WD_E_C}, 
		we have 
	\begin{align}\label{C01_special_WD_formula} 
	v(t,X_t)&=v(0,X_0)+\int_0^t\partial_x v(s,X_s)\,d X^c_s\nonumber\\ 
	&+ \int_{]0,\,t]\times \R} (v(s,X_{s-}+x)-v(s,X_{s-}))\,(\mu^X-\nu^X)(ds\,dx)+A^v(t), 
	\end{align} 
	where 
	$A^v$ is a predictable $(\mathcal F_t)$-orthogonal  process. 
\end{theorem} 

\emph{Proof of Theorem \ref{P_ident}.}
By assumption,  $X$ is  a special weak Dirichlet process satisfying condition \eqref{Sec:WD_CNS}, and  $v$ is a function of class $C^{0,1}$ satisfying the integrability condition \eqref{Sec:WD_E_C}. 
So we are  in the condition to apply 
Theorem \ref{T_C1_dec_specialweak_Dir} 
to $v(t,\,X_t)$. 
We set 
$\varphi(s,x):=v(s,X_{s-}+x)-v(s,X_{s-})$. 
Since $X$ is of finite quadratic variation and verifies \eqref{Sec:WD_CNS}, and $X$ and $v$ satisfy {\eqref{Sec:WD_E_C}, by 
\eqref{A2_Aloc_NEW} and \eqref{F2_Aloc_NEW}
 we see that the process $\varphi$ verifies condition \eqref{cond_int_A} with $A=\{|x| >1\}$. 
Moreover $\varphi(s,0)=0$.   Since  $\mu$ verifies Hypothesis \ref{H_nu} and  $X$ verifies Hypothesis \ref{H_X_mu}, we can apply Proposition \ref{L_ident_mu_muX} to $\varphi(s,x)$. Identity \eqref{C01_special_WD_formula} becomes 
\begin{align} 
v(t,\,X_t) 
&=v(0,X_0) + \int_{]0,\,t]\times \R} (v(s,X_{s-}+ \tilde \gamma(s,e))-v(s,X_{s-}))\,(\mu-\nu)(ds\,de) 
\nonumber\\ 
&\,\,+ \int_{]0,\,t]}\partial_x v(s,X_{s})\,d X^{c}_s+ A^v(t).\label{dec_Y} 
\end{align}

At this point we recall that the process $Y_t = v(t,X_t)$ fulfills the BSDE \eqref{GeneralBSDE}, 
which can be rewritten as 
\begin{align}\label{GeneralBSDE_forward} 
Y_t &= Y_0 + \int_{]0,\,t]} Z_s \, dM_s + \int_{]0,\,t]\times \R} U_s(e)\,(\mu-\nu)(ds\,de)\nonumber\\ 
&\quad - \int_{]0,\,t]} \tilde g(s,\,Y_{s-},\,Z_s)\,  d \zeta_s - \int_{]0,\,t]\times \R} \tilde f(s,\,e,\,Y_{s-},\,U_{s}(e))\,  \lambda(ds\,de). 
\end{align} 
By Proposition 
\ref{Sec:WD_P_unique_decomp}
the uniqueness of the  decomposition \eqref{dec_Y} yields   identity \eqref{id_3} and 
\begin{equation}\label{Id2} 
\int_{]0,\,t]} Z_s \, dM_s=\int_{]0,\,t]} \partial_x v(s,X_{s})\,d X_s^c. 
\end{equation} 
In particular, from \eqref{Id2} we get 
\begin{align*} 
0 &=\langle \int_{]0,\,t]} Z_s dM_s-\int_{]0,\,t]} \partial_x v(s,X_{s})\,d X_s^c,\, M_t \rangle\\ 
&=\int_{]0,\,t]} Z_s d\langle M \rangle_s -\int_{]0,\,t]} \partial_x v(s,X_{s})\,\frac{d \langle X^c,\, M \rangle_s}{d \langle M\rangle_s}\,d \langle M\rangle_s\\ 
&=\int_{]0,\,t]} \bigg(Z_s  - \partial_x v(s,X_{s})\,\frac{d \langle  X^c,\, M \rangle_s}{d \langle M\rangle_s}\bigg)\,d \langle M\rangle_s, 
\end{align*} 
that gives identification \eqref{id_2BIS}. 
 
If in addition we assume that $H\in \mathcal G^2_{\rm loc}(\mu)$, the predictable bracket at time $t$ of the purely discontinuous martingale in   identity \eqref{id_3} is well-defined, and equals $C(H)$ by Theorem 11.21, point 3),   in \cite{chineseBook}. 
Since $C(H)_T=0$ a.s., the conclusion follows from Proposition
\ref{P_forBSDEs2}.
\qed

Let us now consider a BSDE  driven only by a purely discontinuous martingale,  of the form \eqref{GeneralBSDE_disc}. 
We formulate the following assumption for a couple of adapted processes $(X,Y)$. 

We first introduce some notations. Let $E$ be a closed subset of $\R$ on which $X$ takes values. Given a  c\`adl\`ag function $\varphi: [0,\,T] \rightarrow \R$, we denote by $\mathcal C_{\varphi}$ the set of times $t \in [0,\,T]$ for which there is a left (resp. right) neighborhood $I_{t-} = ]t-\varepsilon,\,t[$ (resp. $I_{t+} = [t,\,t+\varepsilon[$) such that $\varphi$ is  constant on $I_{t-}$ and $I_{t+}$.
\begin{hypothesis}\label{H_chain_rule_C0} 
	\hspace{2em} 
	\begin{itemize} 
		\item[(i)] 
		$Y$ is an $(\mathcal F_t)$-orthogonal process such that $\sum_{s \leq T}|\Delta Y_s| < \infty$, a.s.
		\item[(ii)]$X$ is a c\`adl\`ag process and  $Y_t = v(t,\,X_t)$ for 
		some 
		deterministic function $v: [0,T] \times \R \rightarrow \R$, satisfying  the  integrability condition 
		\begin{equation}\label{EC_C0} 
		\int_{]0,\,\cdot]\times \R}|v(t,X_{t-}+x)-v(t,X_{t-})|\,\mu^X(dt\,dx)\,\,\in \mathcal{A}^+_{\rm loc}. 
		\end{equation} 
		\item[(iii)] There exists $\mathcal C \in [0,\,T]$ such that for $\omega$ a.s. $\mathcal C \supset \mathcal C_X(\omega)$, and 
		\begin{itemize}
		\item $\forall t \in \mathcal C$, $t \mapsto v(t,x)$ is continuous $\forall x \in E$;
		\item $\forall t \in \mathcal C^c$, $x \in E$,  $(t, x)$ is a continuity point of $v$.
		\end{itemize}
	\end{itemize} 
\end{hypothesis}
\begin{remark}\label{R_1BIS}
	Item (iii) of Hypothesis \ref{H_chain_rule_C0}  is fulfilled in  two typical situations.
	\begin{enumerate}
	\item$\mathcal C=[0,\,T]$. Almost surely $X$ admits a finite number of jumps and $t \mapsto v(t,x)$ is continuous $\forall x \in E$.
	\item$\mathcal C=\emptyset$ and  $v|_{[0,\,T] \times E}$ is continuous.
	\end{enumerate}
\end{remark}
\begin{remark}\label{R_2BIS}
Assume that Hypothesis \ref{H_chain_rule_C0}-(iii)  holds. Then
\begin{itemize}
\item[(i)]	
$Y_t = v(t, X_t)$ is necessarily a  c\`adl\`ag process.
\item[(ii)] $\forall t \in [0,\,T]$, $\Delta Y_t = v(t,X_t)- v(t, X_{t-})$.
\end{itemize}
\end{remark}

Our second main result is the following. 
\begin{theorem}\label{P_ident_C0} 
	Let  $\mu$ satisfy Hypothesis \ref{H_nu}, and assume that
	  $X$ is a process such that $(X,\mu)$ verifies Hypothesis \ref{H_X_mu}. 
	Let   $(Y,U)$ be a solution to the BSDE \eqref{GeneralBSDE_disc}, 
	such that $(X,Y)$ satisfies Hypothesis  \ref{H_chain_rule_C0} with corresponding function $v$. 
	Then,  the random field  $U$ 	satisfies 
	\begin{equation}\label{id_3disc} 
	\int_{]0,\,t]\times \R} H_s(e)\,(\mu-\nu)(ds\,de)=0 \quad \forall t \in ]0,\,T],\,\,\textup{a.s.}, 
	\end{equation} 
	with 
	\begin{equation} 
	H_s(e):=U_s(e)-(v(s,X_{s-}+ \tilde\gamma(s,e))-v(s,X_{s-})).\label{K} 
	\end{equation} 
	If, in addition, $H \in \mathcal G^2_{\rm loc}(\mu)$, 
	then there exists a predictable process $(l_s)$ such that 
	$$
	H_s(e)= l_s\,\one_{K}(s),\quad d\P\,\nu(ds\,de)\textup{-a.e.}
	$$ 
	In particular, 
\begin{equation}\label{ID_cont_2}
H_s(e)= 0, \quad d\P\,\nu^c(ds\,de)\textup{-a.e.}
\end{equation}
and	
\begin{equation}\label{ID_disc_2}
H_s(e)= l_s\,\one_{K}(s), \quad d\P\,\nu^d(ds\,de)\textup{-a.e.}
\end{equation}
\end{theorem} 

The proof of Theorem \ref{P_ident_C0}  is based on the following stability result for 
c\`adl\`ag processes,  which was  the object  of Proposition 5.31  
in \cite{BandiniRusso1}.

	 
\begin{proposition}\label{P_C00_chain_rule} 
Let $(X,Y)$ be a couple of $(\mathcal F_t)$-adapted processes satisfying Hypothesis  \ref{H_chain_rule_C0} with corresponding function $v$.
Then $v(t,X_t)$ is an $(\mathcal F_t)$-special weak Dirichlet process with decomposition 
\begin{equation}\label{C00_decomp} 
v(t,X_t)=v(0,X_0) 
+ \int_{]0,\,t]\times \R} (v(s,X_{s-}+x)-v(s,X_{s-}))\,(\mu^X-\nu^X)(ds\,dx)+A^v(t), 
\end{equation} 
where $A^v$ is a predictable $(\mathcal F_t)$-orthogonal  process. 
 
\end{proposition}

\emph{Proof of Theorem \ref{P_ident_C0}.}
By assumption, the couple $(X,Y)$ satisfies Hypothesis \ref{H_chain_rule_C0} with corresponding function $v$. 
We are then in the condition to apply 
Proposition \ref{P_C00_chain_rule} 
to $v(t,\,X_t)$, which gives \eqref{C00_decomp}. 
Set 
$ 
\varphi(s,x):=v(s,X_{s-}+x)-v(s,X_{s-})$. 
By condition (ii) in Hypothesis \ref{H_chain_rule_C0}, the process $\varphi$ verifies condition \eqref{cond_int_A} with $A=\Omega \times [0,\,T] \times \R$. 
Moreover $\varphi(s,0)=0$.   Since  $\mu$ verifies Hypothesis \ref{H_nu}, 
and  $(X, \mu)$ verifies Hypothesis \ref{H_X_mu}, 
we can apply Proposition \ref{L_ident_mu_muX} to $\varphi(s,x)$. Identity \eqref{C00_decomp}  becomes 
\begin{align} 
v(t,\,X_t) 
=v(0,X_0) + \int_{]0,\,t]\times \R} (v(s,X_{s-}+ \tilde \gamma(s,e))-v(s,X_{s-}))\,(\mu-\nu)(ds\,de) 
+ A^v(t).\label{dec_YC0} 
\end{align} 
At this point we recall that the process $Y_t = v(t,X_t)$ fulfills BSDE \eqref{GeneralBSDE_disc}, 
which can be rewritten as 
\begin{align}\label{GeneralBSDE_forwardC0} 
Y_t = Y_0 + \int_{]0,\,t]\times \R} U_s(e)\,(\mu-\nu)(ds\,de) - \int_{]0,\,t]\times \R} \tilde f(s,\,e,\,Y_{s-},\,U_{s}(e))\,  \lambda(ds\,de). 
\end{align} 
The uniqueness of  of the canonical decomposition \eqref{dec_YC0} yields   identity 
\eqref{id_3disc}. 
If in addition we assume that $H\in \mathcal G^2_{\rm loc}(\mu)$, the predictable bracket at time $t$ of the purely discontinuous martingale in   identity \eqref{id_3disc} is well-defined, and equals $C(H)$ by Theorem 11.21, point 3),   in \cite{chineseBook}. 
Since $C(H)_T=0$ a.s.,
the conclusion follows from Proposition
\ref{P_forBSDEs2}.
\qed

 \section{Applications}\label{Sec_Appl}

\subsection{BSDEs driven by a jump-diffusion processes}\label{Sec_Ex1} 
	Let us focus on 
	the BSDE  
	\begin{align}\label{BSDE_BaBuPa} 
	Y_t &= g(X_T) + \int_{]t,\,T]} f(s,\,X_s,\,Y_s,\,Z_s,\,U_{s}(\cdot))\,  d s- \int_{]t,\,T]} Z_s \, d W_s - \int_{]t,\,T]\times \R} U_s(e)\,(\mu-\nu)(ds\,de), 
	\end{align} 
	which constitutes	a particular case of the BSDE  \eqref{GeneralBSDE}. 
	This is considered for instance in 
	\cite{BaBuPa}. 
	Here $W$ is a Brownian motion and $\mu(ds\,de)$ is a Poisson random measure with compensator 
	\begin{equation}\label{nu_BBP} 
	\nu(ds \,de)=\lambda(de)\,ds, 
	\end{equation} 
	where  $\lambda$ is a  Borel $\sigma$-finite measure on $\R\setminus \{0\}$  and 
	\begin{equation}\label{integrability_nu} 
	\int_{\R}(1 \wedge |e|^2)\,\lambda(de) < +\infty. 
	\end{equation} 
	Poisson random measures have been introduced for instance in Chapter II, Section 4.b in \cite{JacodBook}. 
	The process $X$ appearing in \eqref{BSDE_BaBuPa}   is a Markov process satisfying the  SDE 
	\begin{equation}\label{EDS_BBP} 
	dX_s = b(X_s)\,ds + \sigma(X_s)\,dW_s + \int_{\R}\gamma(X_{s-},e)\,(\mu-\nu)(ds\,de), \quad s \in [t,\,T], 
	\end{equation} 
	where $b:\R \rightarrow \R$, $\sigma: \R \rightarrow \R$ are globally Lipschitz, and  $\gamma:\R \times \R \rightarrow \R$ is a measurable function such that, for some real $K$, and for all $e \in \R$, 
	\begin{equation}\label{gamma_conditions_BBP} 
	\begin{cases} 
	|\gamma(x,e)|\leq K \,(1\wedge |e|),\quad x \in \R,\\ 
	|\gamma(x_1,e)-\gamma(x_2,e)|\leq K \,|x_1-x_2|\,(1\wedge |e|)\quad x_1,\,x_2 \in \R. 
	\end{cases} 
	\end{equation} 
	For every starting point $x\in \R$ and initial time $t \in [0,\,T]$, there is  a unique solution to \eqref{EDS_BBP} denoted  $X^{t,x}$ (see \cite{BaBuPa}, Section 1). 
	Moreover, modulo suitable assumptions on  the coefficients $(g,f)$,  it is proved  that the BSDE \eqref{BSDE_BaBuPa} admits a unique solution $(Y,Z,U) \in 	 \mathcal{S}^2 \times \mathcal{L}^2\times \mathcal{L}^2(\mu)$,  see Theorem 2.1 in \cite{BaBuPa}, where 
	\begin{align*} 
	\mathcal{S}^2:&=\Big\{(Y_t)_{t \in [0,\,T]}\,\,\textup{adapted c\`adl\`ag :}\,\,\Big \vert \Big \vert\!\sup_{t \in [0,\,T]}|Y_t|\,\Big\vert\Big\vert_{L^2(\Omega)} < \infty\Big\},\\ 
	\mathcal{L}^2:&=\Big\{(Z_t)_{t \in [0,\,T]}\,\,\textup{predictable  :}\,\,\E\Big[\int_0^T Z^2_s \,ds\Big]  < \infty\Big\},
	\end{align*} 
	and $\mathcal{L}^2(\mu)$ is  the space  introduced in \eqref{L2mu}. 
	When $X=X^{t,x}$, the solution $(Y,Z)$ of \eqref{BSDE_BaBuPa} is denoted $(Y^{t,x},Z^{t,x})$. 
	In \cite{BaBuPa} it is proved that 
	\begin{equation}\label{u} 
	v(t,x):=Y_t^{t,x}, \quad (t,x)\in [0,\,T]\times \R, 
	\end{equation} 
	satisfies 
	$ 
	Y_s^{t,x}=v(s,X_s^{t,x}) 
	$ 
	for every $(t,x) \in [0,\,T] \times \R$, $s \in [t,\,T]$. 
 
\begin{lemma}\label{L_BBP_X} 
	Let  $\mu$ and $X$ be respectively the Poisson  random measure and  the stochastic  process satisfying the SDE \eqref{EDS_BBP}. 
	Then $J=K = \emptyset$, $\mu$ satisfies Hypothesis \ref{H_nu} and $(X, \mu)$ 	fulfills Hypothesis \ref{H_X_mu} with  
	\begin{align} 
	X^i_t &= \int_{]0,\,t]\times \R}\gamma(X_{s-},e)\,(\mu-\nu)(ds\,de),\label{Xi_BBP}\\ 
	X^p_t &= \int_{0}^{t}b(X_s)\,ds + \int_{0}^{t}\sigma(X_s)\,dW_s,\label{Xp_BBP} \\
	\tilde \gamma(\omega,s,e)&=\gamma(X_{s-}(\omega),e).\nonumber 
	\end{align} 
\end{lemma} 
\proof 
Our aim is to apply Lemma \ref{Ex_guida}. 
We start by noticing that $\nu$ in \eqref{nu_BBP} is in the form \eqref{nu_dis} with $A_s=s$.
Therefore by Remark \ref{R_H_mu_BIS} item 2. and successively item 1, $J=K=\emptyset$, and Hypothesis \ref{H_nu} is  verified. 
On the other hand, the process $X$ satisfies the stochastic differential 
equation \eqref{EDS_BBP}, which is a particular case of 
\eqref{X_SDE}  when $B_s=s$, $N_s=W_s$, and  $b$, $\sigma$, $\gamma$ are  time homogeneous. 
$b$ and $\sigma$ verify 
\eqref{Ver_1}, \eqref{Ver_2} since they have linear growth. 
Condition  \eqref{Ver_3} can be verified using the characterization of $\mathcal G^1_{\rm loc}(\mu)$ in 
Theorem 1.33, point c), Chapter II,  in \cite{JacodBook}:  for a given predictable random field $W$ defined on $\tilde \Omega$ such that $\hat W=0$,  that theorem specifies that, whenever $|W|^2 \one_{\{|W|\leq 1\}} \star \nu+|W| \one_{\{|W|> 1\}} \star \nu \in \mathcal A_{\rm loc}^+$, then $W \in\mathcal G^1_{\rm loc}(\mu)$.
That property  follows from 
\eqref{integrability_nu} and \eqref{gamma_conditions_BBP}. 
 
Then, by Lemma \ref{Ex_guida},  $X$ verifies Hypothesis   
\ref{H_X_mu} with decomposition $X=X^i+X^p$,  where $X^i$ and $X^p$ are given respectively by  \eqref{Xi_BBP} and \eqref{Xp_BBP}, and with 
$\tilde \gamma(\omega, s,e) = \gamma(X_{s-}(\omega),e)$.
\endproof 

We aim at applying
 Theorem \ref{P_ident_C0} to 
 BSDE \eqref{BSDE_BaBuPa}.
To this end, we need the following preliminary result. 
\begin{lemma}\label{L_BBP_EC} 
	Let  $\mu$ and $X$ be respectively the Poisson  random measure and  the stochastic  process satisfying the SDE \eqref{EDS_BBP}. 
	Let 
	$v:[0,\,T] \times \R \rightarrow \R$ be a function of  $C^{0,1}$ class such that $x \mapsto \partial_x v(s,x)$ has linear growth, uniformly in $s$. 
	Then condition \eqref{Sec:WD_E_C} holds for $X$ and $v$. 
\end{lemma} 
\proof 
We have 
\begin{align} 
&\int_{]0,\cdot]\times \R} |v(s,X_{s-} +x )-v(s,X_{s-})-x\,\partial_x v(s,X_{s-})|\,\one_{\{|x| >1\}}\,\mu^X(ds\,dx)\nonumber\\ 
&=\sum_{0<s \leq \cdot} |v(s,X_{s})-v(s,X_{s-})-\partial_x v(s,X_{s-})\,\Delta X_s|\,\one_{\{|\Delta X_s| >1\}}\nonumber\\ 
&\leq \sum_{0<s \leq \cdot}|\Delta X_s|\,\one_{\{|\Delta X_s| >1\}}\,\left(\int_0^1|\partial_{x} v(s,X_{s-}+ a\,\Delta X_s)|\,da +\int_0^1 |\partial_{x} v(s,X_{s-})|\,da\right)\nonumber \\ 
&\leq 2\,C\,\sum_{0<s \leq \cdot}\,|X_{s-}||\Delta X_s|\,\one_{\{|\Delta X_s| >1\}} + \sum_{s \leq t}|\Delta X_s|^2\,C\,\one_{\{|\Delta X_s| >1\}} \nonumber\\ 
&= 2\,C\,\int_{]0,\cdot]\times \R}|X_{s-}|\,|x|\,\one_{\{|x| >1\}} \,\mu^X(ds\,dx) + \sum_{s \leq \cdot}|\Delta X_s|^2\,\one_{\{|\Delta X_s| >1\}},\label{verificaAloc} 
\end{align} 
for some constant $C$.
Since $X$ is of finite quadratic variation, by Lemma 3.9-(ii) in \cite{BandiniRusso1} we have that $\sum_{s \in]0,\,T]}|\Delta X_s|^2 < \infty$, a.s.
Consequently, the second term in the right-hand side of \eqref{verificaAloc} belongs to $\mathcal A_{\rm loc}^+$ 
if we prove that
\begin{equation}\label{sum_Aloc} 
\sum_{s \in]0,\,\cdot]}|\Delta X_s|^2\in \mathcal A_{\rm loc}^+.
\end{equation} 
 Since by \eqref{EDS_BBP} $\Delta X_s= \int_{\R}\gamma(X_{s-},e)\,\mu(ds\,de)$, we have 
$$ 
\sum_{s \in]0,\,\cdot]}|\Delta X_s|^2 = \sum_{s \in]0,\,\cdot]}\left|\int_{\R}\gamma(X_{s-},e)\,\mu(ds\,de)\right|^2=\int_{]0,\cdot]\times \R}|\gamma(X_{s-},e)|^2\,\mu(ds\,de), 
$$ 
and \eqref{sum_Aloc} reads 
\begin{equation}\label{sum_Aloc2} 
\int_{]0,\cdot]\times \R}|\gamma(X_{s-},e)|^2\,\mu(ds\,de)\in \mathcal A_{\rm loc}^+. 
\end{equation} 
 The integral in the left-hand side of \eqref{sum_Aloc2}  exists almost surely. Indeed, 
$|\gamma(x,e)|\leq K\,(1 \wedge |e|)$ for every $x \in \R$, 
$ 
\int_{\R}(1 \wedge |e|^2) \,\lambda(de) < \infty 
$ 
(see, respectively,  \eqref{gamma_conditions_BBP} and \eqref{integrability_nu}).
Since it is c\`agl\`ad, then it is 
 locally bounded, see for instance the lines above Theorem 15, Chapter IV, in \cite{protter}. Consequently, it belongs to $\mathcal A_{\rm loc}^+$.

Finally, the first term in the right-hand side of \eqref{verificaAloc} belongs to $\mathcal A^+_{\rm loc}$, 
taking into account \eqref{Sec:WD_CNS} and the fact  that $X_{s-}$ is locally bounded 
being c\`agl\`ad. The conclusion follows. 
\endproof 
We are ready to give the identification result.
\begin{corollary}\label{C_id_BBP} 
	Let   $(Y,Z,U)\in 	\mathcal{S}^2 \times \mathcal{L}^2\times \mathcal{L}^2(\mu)$ be the unique solution  to the BSDE \eqref{BSDE_BaBuPa}. If  the function $v$ defined in \eqref{u} is of class $C^{0,1}$  such that $x \mapsto \partial_x v(t,x)$ has linear growth, uniformly in $t$, 
	then  the pair $(Z,U)$ 
	satisfies 
	\begin{equation}\label{Z_BBP} 
Z_t = \sigma(X_t)\,\partial_x u(t,X_t) \quad d\P\, dt\textup{-a.e.,} 
	\end{equation} 
	\begin{equation}\label{id_3BBP} 
	\int_{]0,\,t]\times \R} H_s(e)\,(\mu-\nu)(ds\,de)=0, \quad \forall t \in ]0,\,T],\,\, \textup{a.s.} 
	\end{equation} 
	where 
	\begin{equation}\label{K_BBP} 
	H_s(e):=U_s(e)-(v(s,X_{s-}+ \gamma(s,X_{s-},e))-v(s,X_{s-})). 
	\end{equation} 
	If in addition $H \in \mathcal G^2_{\rm loc}(\mu)$,  
	\begin{equation} 
	U_s(e) = v(s,X_{s-}+ \gamma(s,X_{s-},e))-v(s,X_{s-})\quad  d\P\,\lambda(de)\,ds \textup{-a.e.}\label{U_BBP} 
	\end{equation}	 
\end{corollary}	 
\proof 
We aim at applying Theorem \ref{P_ident}. 
By Lemma \ref{L_BBP_X}, $\mu$  satisfies Hypothesis \ref{H_nu} 
with $\tilde \gamma(s,e) 
= \gamma(s,X_{s-},e)$. 
Since $X$ is a special semimartingale, then  condition \eqref{Sec:WD_CNS} holds by  Corollary 11.26 in \cite{chineseBook}. Moreover,  $X$ is obviously a special weak Dirichlet process with  finite quadratic variation. 
By Lemma \ref{L_BBP_EC}, condition \eqref{Sec:WD_E_C} holds for $X$ and  $v$, which implies that Hypothesis \ref{H_chain_rule_C01} is verified. 
 
We can then apply Theorem \ref{P_ident}: since $X^c_\cdot=\int_0^\cdot \sigma(X_t)\,dW_t$ and $M=W$, 
\eqref{id_2BIS}    gives \eqref{Z_BBP}, while \eqref{id_3}-\eqref{K_def}  with $\tilde \gamma(s,e)= \gamma(s, X_{s-},e)$ yield \eqref{id_3BBP}-\eqref{K_BBP}. 
If in addition  $H \in \mathcal{G}^2(\mu)$,   since $J=K=\emptyset$,  see Lemma \ref{L_BBP_X}, $\nu=\nu^c$ and 
\eqref{U_BBP} follows by \eqref{ID_cont_1}. 
\endproof 
\begin{remark}\label{R_Fuhrman_Tessitore} 
	When the BSDE  \eqref{BSDE_BaBuPa} is driven only by a standard Brownian motion,  an  identification result for $Z$ analogous to \eqref{Z_BBP} 
	has  been established by \cite{FuhrmanTessitore}, even supposing only that $f$ is  Lipschitz with respect to $Z$. 
\end{remark} 
 
\begin{remark}\label{R_distributional_drift} 
Theorem \ref{P_ident} also potentially  applies to the 
 cases of  BSDEs
driven by a continuous martingale when the underlying process
$X$ is a solution of an SDE with singular (distributional) 
drift. 
In the literature there are plenty of cases of (even continuous) Markov processes that are not semimartingales.  
Typical examples of such underlying processes $X$
 are solutions of an SDE with distributional drift, see e.g. 
	\cite{frw1},  \cite{rtrut},  \cite{issoglio}, of the type 
	\begin{equation}\label{SDE_distr} 
	dX_t = \beta(X_t)\,dt + d W_t, 
	\end{equation} 
	for a class of Schwartz distributions $\beta$. 
	In this case $X$ 
	is generally not a semimartingale but only a Dirichlet process, so that, 
	for $v \in C^{0,1}$, $v(t,X_t)$ is a (special) weak Dirichlet process. 
	Forward BSDEs related to an underlying process $X$ solving \eqref{SDE_distr} 
	have been  studied  for instance in \cite{wurzer}, when the terminal type 
	is random. 
However we do not perform 
a more refined analysis of examples in that direction since it goes  beyond the scope of the paper.
\end{remark}

 
 
 \subsection{On a class of BSDEs driven by a quasi-left-continuous random measure}
	In \cite{CoFu-m} the authors study a   BSDE driven by an integer-valued random measure $\mu$ associated to a given pure jump Markov process $X$, 
	of the   form 
	\begin{align}\label{BSDE_CoFu_jump} 
	Y_t &= g(X_T) + \int_{]t,\,T]} f(s,\,X_s,\,Y_s,\,U_{s}(\cdot))\,  d s - \int_{]t,\,T]\times \R} U_s(e)\,(\mu-\nu)(ds\,de). 
	\end{align} 
	The underlying process $X$ is generated by a marked point process $(T_n, \zeta_n)$, where $(T_n)_n$ are increasing random times such that 
	$ 
	T_n \in ]0,\,\infty[, 
	$ 
	where either there is a finite number of 
 times  $(T_n)_n$  or $\lim_{n \rightarrow \infty} T_n = + \infty$, 
	and $\zeta_n$ are random variables in $\R$, see e.g. Chapter III, Section 2 b., in \cite{jacod_book}. This means that $X$ is a c\`adl\`ag process such that   $X_{t}= \zeta_n$ for $t \in [T_n,\,T_{n+1}[$, for every $n \in \N$. In particular, $X$ has a finite   number of jumps on each compact. 
	The associated  integer-valued random measure   $\mu$ is the sum of the Dirac measures concentrated at the marked point process $(T_n,\zeta_n)$, and can  be written as 
	\begin{equation}\label{4.12_BIS} 
	\mu(ds\,de) = \sum_{s \in [0,\,T]}\one_{\{X_{s-}\neq X_s\}}\,\delta_{(s,X_s)}(dt\,de). 
	\end{equation} 
	Given a measure $\mu$ in the form \eqref{4.12_BIS}, it 
	is related to the jump measure $\mu^X$ in the following way: 
	for every Borel subset $A$ of $\R$, 
	\begin{equation}\label{change_meas} 
	\int_{]0,\,T] \times \R} \one_{A}(e-X_{s-})\,\mu(ds\,de) = \int_{]0,\,T] \times \R} \one_{A}(x)\,\mu^X(ds\,dx). 
	\end{equation} 
	This is for instance explained in  Example 3.22 in \cite{jacod_book}. 
	The pure jump process $X$ then satisfies the equation 
	\begin{align} 
	\label{X_pure_jump} 
	X_t =  X_0 + \sum_{0<s \leq t} \Delta X_{s}=X_0 + \int_{]0,t]\times \R}(e-X_{s-})\,\mu(ds\,de). 
	\end{align} 
	The compensator of $\mu(ds\,de)$ is 
	\begin{equation}\label{nu_CoFu} 
	\nu(ds\,de)=\lambda(s,X_{s-},de)\,ds, 
	\end{equation} 
	where  $\lambda$  the  is the transition rate measure of the process satisfying 
	\begin{equation}\label{bounded_rate} 
	\sup_{t\in [0,\,T],\,x\in \R}\lambda(t,x,\R) <\infty, 
	\end{equation} 
	see Section 2.1 in \cite{CoFu-m}. 
	 
	Under suitable assumptions  on  the coefficients $(g,f)$, Theorem 3.4 in \cite{CoFu-m} states that the BSDE \eqref{BSDE_CoFu_jump} admits a unique solution $(Y,U) \in 	\mathcal{L}^2 \times \mathcal{L}^2(\mu)$,    where 
	$\mathcal{L}^2(\mu)$ and $\mathcal{L}^2$ are the spaces  introduced in Section \ref{Sec_Ex1}.
	Theorem 4.4 in \cite{CoFu-m} shows moreover that there exists a measurable function  
	$v: [0,\,T]\times \R \rightarrow \R$  such that 
	\begin{align} 
	&\forall\,\, 
	e \in \R,\,\,  t \mapsto v(t,e)\,\, \textup{is absolutely continuous on} \,\,[0,\,T],\label{prop_u_1}\\ 
	& (v(s,X_s)) \in \mathcal L^2\,\,\textup{and}\,\,(v(s,e)-v(s,X_{s-}), \,(s,e) \in [0,\,T] \times \R)   \in \mathcal L^2(\mu),
	\label{prop_u_2} 	
	\end{align} 
	and the unique solution of the BSDE  \eqref{BSDE_CoFu_jump} can be represented as 
	\begin{align} 
	Y_s&=v(s,X_s),\quad s\in [0,\,T].\label{uY}
	\end{align} 
\begin{lemma}\label{L_CoFu_X} 
	Let $\mu$ be the integer-valued random measure  in   \eqref{4.12_BIS} with compensator $\nu$ given by \eqref {nu_CoFu},  and  $X$ be  the associated  pure jump Markov process satisfying \eqref{X_pure_jump}. 
	Then $J=K=\emptyset$, $\mu$ satisfies Hypothesis \ref{H_nu} and $(X, \mu)$ fulfills Hypothesis \ref{H_X_mu} 
	with  $\tilde \gamma(\omega, s,e)= e-X_{s-}(\omega)$. 
\end{lemma} 
\proof	 
Since $\nu$ in \eqref{nu_CoFu} is in the form \eqref{nu_dis} with $A_s=s$, 
by Remark \ref{R_H_mu_BIS}, item 2. and successively  item 1., $J=K=\emptyset $ and Hypothesis \ref{H_nu} is  verified.
Let us now prove that $(X, \mu)$ fulfills Hypothesis \ref{H_X_mu} 
	with  $\tilde \gamma(\omega, s,e)= e-X_{s-}(\omega)$.  
We observe that \eqref{4.12_BIS} implies that  $\{\Delta X \neq  0 \}=D$.
Therefore, recalling that by Remark \ref{R_H_mu_BIS}-1., $D= \cup_n [[T^i_n]]$
for some sequence of totally inaccessible times} $(T^i_n)_n$ such that  $[[T^i_n]] \cap [[T^i_m]] = \emptyset$, $n \neq m$,
 we see that $X=X^i$ is quasi-left-continuous.
 
Finally,  by definition of $\mu$ we have 
$$ 
\sper{\int_{]0,\,T]\times \R}\mu(ds\,de)\,| (e-X_{s-})-\Delta X_s|}=0, 
$$ 
therefore $X^i$ satisfies Hypothesis \ref{H_X_mu}-1. with  
$\tilde \gamma(\omega, s,e)= e-X_{s-}(\omega)$. Since $X^p=0$, Hypothesis \ref{H_X_mu}-2. trivially holds. 
\endproof 

Let us now apply Theorem \ref{P_ident_C0}  to the present framework. 
We start with a  preliminary observation. 
\begin{lemma}\label{L_CoFu_EC} 
	Let $\mu$ be the integer-valued random measure  in   \eqref{4.12_BIS} with compensator $\nu$ given by \eqref {nu_CoFu},  and  $X$ be  the associated  pure jump Markov process satisfying \eqref{X_pure_jump}.  
	Let $v:[0,\,T]\times \R \rightarrow \R$ be a 
	function satisfying 
	\eqref{prop_u_2}
	 and such that, for $E = \R$,
	 \begin{equation}\label{new_continuity}
\forall\, 
	e \in E, \,\,  t \mapsto v(t,e)\,\,  \textup{is  continuous on } [0,\,T].
		 \end{equation}
	Suppose that  $Y_t = v(t, X_t)$ is an $(\mathcal F_t)$-orthogonal process.
	Then $(X,Y)$ satisfies Hypothesis \ref{H_chain_rule_C0} with corresponding function $v$. 
\end{lemma} 
\proof	 
By Remark \ref{R_1BIS}, from \eqref{new_continuity} it follows that Hypothesis \ref{H_chain_rule_C0}-(iii) is verified.
Taking into account Remark \ref{R_2BIS},
 it follows that $Y_t=v(t, X_t)$ has a finite number of jumps, 
 so that $\sum_{s \leq T}| \Delta Y_s| < \infty$ a.s. In particular, Hypothesis \ref{H_chain_rule_C0}-(i) holds. 
 
To verify the validity of condition (ii) of Hypothesis \ref{H_chain_rule_C0},
 we have to show that \eqref{EC_C0} holds. 
Denoting $||\lambda||_{\infty}= \sup_{t\in [0,\,T],\,x \in \R}|\lambda(t,x,\R)|$, by \eqref{change_meas} we have 
\begin{align*} 
&\sper{\int_{]0,\,T] \times \R}|v(s,X_{s-}+x)-v(s,\,X_{s-})|\,\mu^X(ds\,dx)}\\ 
&=\sper{\int_{]0,\,T] \times \R}|v(s,e)-v(s,\,X_{s-})|\,\mu(ds\,de)}\\ 
&=\sper{\int_{]0,\,T] \times \R}|v(s,e)-v(s,\,X_{s-})|\,\lambda(s,X_{s-},\,de)\,ds}\\ 
&\leq T\,||\lambda||^{1/2}_{\infty}\,||v(s,e)-v(s,\,X_{s-})||^{1/2}_{\mathcal L^2(\mu)} 
\end{align*} 
and  the conclusion follows since  $v(s,e)-v(s, X_{s-})\in \mathcal L^2(\mu)$ by \eqref{prop_u_2}. 
\endproof 
We have the following identification result.
\begin{corollary}\label{C_id_CoFu} 
	Let   $(Y,U)\in 	\mathcal{L}^2 \times  \mathcal{L}^2(\mu)$ be the unique solution  to the BSDE \eqref{BSDE_CoFu_jump} and  $X$, $v$  be respectively the process and the function  appearing  in \eqref{uY}. 
	Then 
	the random field  $U$	satisfies 
	\begin{equation} 
	U_t(e) = v(t,e)-v(t,X_{t-})\quad  d\P\,\lambda(t,X_{t-},\,de)\,dt \textup{-a.e.}\label{U_CoFu} 
	\end{equation}	 
\end{corollary} 
\proof 
We aim at applying  Theorem \ref{P_ident_C0}. 
By Lemma \ref{L_CoFu_X}, $\mu$  satisfies Hypothesis \ref{H_nu} and $(X, \mu)$ fulfills   Hypothesis \ref{H_X_mu} 
with $\tilde \gamma(s,e) 
= e - X_{s-}$. 
Moreover, 
by Lemma \ref{L_CoFu_EC} and Remark \ref{R_orth_BSDE}, $(X,Y)$ satisfies   Hypothesis \ref{H_chain_rule_C0} with corresponding function $v$. 
We can then apply Theorem \ref{P_ident_C0}. 
We have 
\begin{align}\label{K_CoFu} 
H_s(e)&:=U_s(e)-(v(s,X_{s-}+ \tilde \gamma(s,e))-v(s,X_{s-}))\nonumber\\ 
&=U_s(e)-(v(s,e)-v(s,X_{s-})), 
\end{align} 
which belongs to $\mathcal L^2(\mu)$, and therefore to $\mathcal{G}^2(\mu)$. Recalling that  $J = K= \emptyset$,   see Lemma \ref{L_CoFu_X}, $\nu=\nu^c$ and
\eqref{U_CoFu} follows by 
\eqref{ID_cont_2}. 
\endproof 	 

\begin{remark}
The result in Corollary \ref{C_id_CoFu}  is not new, since it  retrieves with a different method, without using the specific form of $v$,  the result obtained in \cite{CoFu-m}. In particular, our identification does  not need to use the absolute continuity property \eqref{prop_u_1}.
\end{remark}

\subsection{On a class of BSDEs driven by a non quasi-left-continuous random measures} 
 
In the recent paper   \cite{BandiniBSDE}, the author studies the existence and uniqueness for a BSDE driven by a purely discontinuous martingale of the form 
\begin{align}\label{BSDE_Ba} 
Y_t &= \xi + \int_{]t,\,T]} \tilde f(s,\,Y_{s-},\,U_{s}(\cdot))\,  d A_s - \int_{]t,\,T]\times \R} U_s(e)\,(\mu-\nu)(ds\,de),
\end{align} 
for given  data $\xi$, $\tilde f$.
Here $\mu(ds\,de)$ is  an integer-valued random measure with  compensator $\nu(ds\,de)=d A_s \,\phi_s(de)$, where  $\phi$ is a probability kernel and  $A$ is a right-continuous nondecreasing predictable process, such that $\hat{\nu}_s(\R)=\Delta {A}_s \leq 1$ for every $s$. 
For any positive constant $\beta$,  $\mathcal E^\beta$ will denote the Dol\'eans-Dade exponential of the process $\beta A$. 
We consider the  weighted spaces 
\begin{align*} 
&\mathcal{L}^2_{\beta}(A):=\Big\{\textup{adapted c\`adl\`ag  processes}\,\,(Y_s)_{s \in [0,\,T]},\,\,\textup{s.t.}\,\E\big[\int_0^T \mathcal E^\beta_s |Y_{s-}|^2 \,dA_s\big]  < \infty\Big\},\\ 
&\mathcal{G}^2_{\beta}(\mu):=\Big\{\textup{predictable  processes}\,\,(U_s(\cdot))_{s \in [0,\,T]},\,\,\textup{s.t.}\\ 
&||U||^2_{\mathcal{G}^2_{\beta}(\mu)}:=\E\big[\int_{]0,\,T]\times \R}\mathcal E^\beta_s\,  |U_s(e)-\hat{U}_s|^2 \,\nu(ds\,de) + \sum_{s \in ]0,\,T]}\mathcal E^{\beta}_s\,|\hat{U}_s|^2(1-\Delta A_s)\big]  < \infty\Big\}. 
\end{align*} 
In \cite{BandiniBSDE} the author  considers solutions
$(Y,U) \in \mathcal{L}^2_{\beta}(A) \times \mathcal{G}^2_{\beta}(\mu)$. 
Suitable assumptions  are required on the triplet $(\tilde f,\xi, \beta)$. 
In particular $\tilde f$ is of Lipschitz type in the third and fourth variable and $\xi$ is a square integrable random variable with some weight. 
Moreover,  the following  technical assumption has to be fulfilled: there exists $\varepsilon \in ]0,\,1[$ such that 
\begin{equation}\label{Cond_BSDE_randmeas} 
2\,|L_y|^2\,|\Delta A_t|^2 \leq 1-\varepsilon, \quad \P{\rm -a.s.},\quad \forall \,t \in [0,\,T], 
\end{equation} 
where $L_y$ is the Lipschitz constant of $\tilde f$ with respect to $y$. 
Under these hypotheses, for $\beta$ large enough, it can be proved that there exists a (unique) solution $(Y,U) \in \mathcal{L}^2_{\beta}(A) \times \mathcal{G}^2_{\beta}(\mu)$ to BSDE \eqref{BSDE_Ba}, see Theorem 4.1 in  \cite{BandiniBSDE}. 
In that theorem, one shows that,
given two solutions $(Y,U)$, $(Y',U')$, then we have $Y_t=Y_t'$ $d\P \,d A_t$-a.e. and $||U-U'||^2_{\mathcal{G}^2_{\beta}(\mu)}=0$. This implies  that $||U-U'||^2_{\mathcal{G}^2(\mu)}=0$, and so, by Remark \ref{R_uniq_G2}, there is a predictable process $(l_t)$ such that 
$U_t(e)-U'_t(e) = l_t\,\one_K(t)$, $d\P\,\nu(dt\,de)$-a.e.
\paragraph{The PDMPs case.}
	Let us now consider a particular case of BSDE \eqref{BSDE_Ba}, namely a BSDE driven by the integer-valued random measure $\mu$ associated to a given  Markov process $X$, 
	of the   form 
	\begin{equation}\label{BSDE_PDP} 
	Y_t = g(X_T) + \int_{]t,\,T]} f(s,\,X_{s-},\,Y_{s-},\,U_{s}(\cdot))\,  d A_s - \int_{]t,\,T]\times \R} U_s(e)\,(\mu-\nu)(ds\,de). 
	\end{equation} 
	We assume that $X$ is a  piecewise deterministic Markov process (PDMP)  associated to a random measure $\mu$,  with values in the interval $]0,1[$.
The   process $X$ is generated by a marked point process $(T_n, \zeta_n)$, where $(T_n)_n$ are increasing random times such that 
	$ 
	T_n \in ]0,\,\infty[, 
	$ 
	where either there is a finite number of 
 times  $(T_n)_n$  or $\lim_{n \rightarrow \infty} T_n = + \infty$, 
	and $\zeta_n$ are random variables in $]0,1[$.

	We will follow the notations in \cite{Da-bo}, Chapter 2, Sections 24 and 26. 	The behavior of the PDMP $X$ is described by a triplet of 
	local characteristics  $(h,\lambda,P)$: 
	$h: ]0,\,1[ \rightarrow \R$  is a Lipschitz continuous function, 
	$\lambda: ]0,1[ \rightarrow \R$ is a measurable function satisfying 
\begin{equation}\label{bounded_lambda} 
	\sup_{x \in ]0,1[}|\lambda(x)| < \infty, 
	\end{equation}
	and $P$ is a transition probability  measure on $[0,1]\times\mathcal{B}(]0,1[)$. 	Some other technical assumptions  are specified in the over-mentioned reference, that we not recall here.
	Let us denote by $\Phi(s,x)$  the unique solution of $
	{g}'(s) = h(g(s))$, $g(0)= x$.
	The process $X$ can  be   defined as 
		\begin{equation} \label{X_eq}
			X(t)=
	\left\{
	\begin{array}{ll}
	\Phi(t,x),\quad t \in [0,\,T_{1}[\\
	\Phi(t-T_n, \zeta_n),\quad t \in [T_n,\,T_{n+1}[.
	\end{array}
	\right.
	\end{equation}
	Set  $N_t = \sum_{n\in\N} \one_{t \geq T_n}$. 
	By Proposition 24.6 in \cite{Da-bo}, 
	we have 
	\begin{equation}\label{Nt_integr} 
	\sper{N_t}< \infty\quad \forall \,t \in \R_+. 
	\end{equation} 
	Notice that the PDMP $X$ verifies the equation 
	\begin{align} 
	X_t=X_0 + \int_{0}^t h(X_s)\,ds + \sum_{0<s \leq t}\Delta X_s.\label{PDP_dynamic} 
	\end{align} 
	In particular $X$ admits a finite number of jumps on each compact interval. 
	By  (26.9) in \cite{Da-bo}, the  random measure   $\mu$  is 
	\begin{align}	 
	\mu(ds\,de) 
	=\sum_n \one_{\{\zeta_{n}\in ]0,1[\}} \delta_{( T_n,\,\zeta_n)}(ds\,de)
	=\sum_{0<s \leq t} 
	\one_{\{X_{s-}\neq X_{s}\}}\delta_{( s,\,X_s)}(ds\,de),\label{mu_PDP} 
	\end{align} 
	which is of the type of \eqref{4.12_BIS}. 
	This implies the validity of \eqref{change_meas}, so that 
	\eqref{PDP_dynamic} can be rewritten as 
	\begin{align} 
	X_t=X_0 + \int_{0}^t h(X_s)\,ds + \int_{]0,\,t]\times ]0,1[}(e-X_{s-})\,\mu(ds\,de).\nonumber 
	\end{align} 
	In the following, by abuse of notations, 
	$\mu$ will denote the trivial extension of previous measure  to the real line in the space variable $e$. 
	In particular 
	\eqref{PDP_dynamic} can be reexpressed as 
	\begin{equation}\label{PDP_dynamic_estended} 
	X_t=X_0 + \int_{0}^t h(X_s)\,ds + \int_{]0,\,t]\times \R}(e-X_{s-})\,\mu(ds\,de). 
	\end{equation} 
	 
	The knowledge of $(h,\,\lambda,\,P)$ completely specifies the dynamics of $X$, see Section 24 in \cite{Da-bo}. 
	According to (26.2) in \cite{Da-bo}, the compensator of $\mu$ has the form 
	\begin{equation}\label{nuPDPs} 
	\nu(ds\, de) = (\lambda(X_{s-})\,ds + d p^{\ast}_s)\,P(X_{s-},\,de), 
	\end{equation} 
	where $\lambda$ has been trivially extended to $[0,1]$ by the zero value, and  
	\begin{equation}\label{p_ast} 
	p^{\ast}_t = \sum_{n =1}^{\infty} \one_{\{t \geq T_n\}}\,\one_{\{X_{{T_n}-} \in \{0,1\}\}} 
	\end{equation} 
	is the predictable process counting the number of jumps of  $X$ from the  boundary of its domain. 
	 
	From \eqref{nuPDPs} we can choose $A_s$ and $\phi_s(de)$ such that 
	$d A_s = \lambda(X_{s-})\,ds + d p^{\ast}_s$ and $\phi_{s}(de)= P(X_{s-},de)$. 
	In particular, $A$ is predictable (not deterministic) and discontinuous, with jumps 
	\begin{equation}\label{Delta_A} 
	\Delta A_s(\omega) =\hat \nu_s(\omega,\R)= \Delta p^{\ast}_s(\omega)=\one_{\{X_{s-}(\omega) \in \{0,1\}\}}. 
	\end{equation} 
	Consequently, $\hat \nu_t(\omega,\R) >0$ if and only if $\hat \nu_t(\omega,\R) =1$, so that 
	\begin{equation}\label{PDP_J=K} 
	J=\{(\omega,t): \hat \nu_t(\omega,\R) >0\}=\{(\omega,t): \hat \nu_t(\omega,\R) =1\}=K, 
	\end{equation} 
	and 
	\begin{equation}\label{K_PDP} 
	K= \{(\omega,t):X_{t-}(\omega)\in \{0,1\}\}. 
	\end{equation} 

\begin{lemma}\label{L_Ban_int} 
	Let $X$ be  the PDMP process with local characteristics $(h,\lambda, P)$,  satisfying \eqref{PDP_dynamic}. 
	Then 
	$$ 
	\int_{]0,\,\cdot]\times \R}|e-X_{s-}|\,\nu(ds\,de) \in \mathcal{A}^+_{\rm loc}. 
	$$ 
\end{lemma} 
\proof 
We start by noticing that 
$\int_{]0,\,T]\times \R}|e-X_{s-}|\,\nu(ds\,de) < \infty$, a.s.
Indeed 
\begin{align*} 
\int_{]0,\,T]\times \R}|e-X_{s-}|\,\nu(ds\,de) &= \int_{]0,\,T]\times ]0,\,1[}|e-X_{s-}|\,(\lambda(X_{s-})\,ds+ d p^{\ast}_s)\,P(X_{s-},de)\\ 
&\leq  
||\lambda||_{\infty}\,(T+ p^{\ast}_T). 
\end{align*} 
For every $t\in [0,\,T]$, the jumps of the process 
$
\Gamma_t:=\int_{]0,\,t]\times \R}|e-X_{s-}|\,\nu(ds\,de) 
$ 
are given by 
$$ 
\Delta \Gamma_t:=\int_{]0,\,1[}|e-X_{t-}|\,\hat \nu_t(de) \leq \hat \nu_t(\R) \leq 1. 
$$ 
Since  $\Gamma_t$ has bounded jumps, it is a locally bounded process  and therefore it belongs to $\mathcal A_{\rm loc}^+$, see for instance the proof of Corollary at page 373 in \cite{protter}. 
\endproof 
 
\begin{lemma}\label{L_Bandini_X} 
	Let  $\mu$ and  $X$ be respectively the  random measure  and  the associated  PDMP with local characteristics $(h,\lambda,P)$ satisfying equation \eqref{PDP_dynamic_estended}. 
	Assume in addition that there exists a  function $\beta: \{0,1\}\rightarrow ]0,1[$, such that 
	\begin{equation}\label{cond_PDP} 
	P(x,\,de)=\delta_{\beta(x)}(de)\quad \textup{a.s.} 
	\end{equation} 
	Then $\mu$ satisfies Hypothesis \ref{H_nu} and $(X, \mu)$ 	fulfills Hypothesis \ref{H_X_mu} with 
	\begin{align} 
	X^i_t &= \int_{]0,\,t]\times \R}(e - X_{s-})\,(\mu-\nu)(ds\,de),\label{Xi_PDP}\\ 
	X^p_t &= X_0 + \int_{0}^{t}h(X_s)\,ds+ \int_{]0,\,t]}\left(\int_{\R}(e - X_{s-})\,P(X_{s-},\,de)\right)(\lambda(X_{s-})\,ds+dp^\ast_s),\label{Xp_PDP} \\
	\tilde \gamma(\omega,s,e)&=(e-X_{s-}(\omega))\,\one_{\{ X_{s-}(\omega)\in ]0,1[\}}(\omega,s). \nonumber
	\end{align} 
\end{lemma} 
\begin{remark}
Condition \eqref{cond_PDP}  implies that 	
	\begin{equation}\label{cond_PDP_BIS} 
	X_s = \beta(X_{s-})\quad \textup{on}\,\,\{(\omega,s): X_{s-}(\omega)\in \{0,1\}\}.
	\end{equation}
\end{remark}

\proof	 
Let us prove that Hypothesis \ref{H_nu}-(i) holds. 
We recall that the measure $\mu$ was characterized by \eqref{mu_PDP}. 
Moreover, given 
$\mu^c= \mu\,\one_{J^c}$ and 
$\nu^c= \nu\,\one_{J^c}$, 
$\nu^c$ is the compensator of    $\mu^c$, see paragraph b) in \cite{jacod_repr}. Taking into account \eqref{nuPDPs}, \eqref{Delta_A} and  \eqref{PDP_J=K}, we have 
\begin{align} 
\nu^{c}(ds\,de) &= \lambda(X_{s-})\,P(X_{s-},de)\,ds.\label{nuc_PDP} 
\end{align} 
We remark that $D\,\cap J^c = \{(\omega, t) : \, \mu^c(\omega, \{t\} \times \R)>0\}$.
By Remark \ref{R_H_mu}-(ii),  we have
$D\,\cap J^c= \cup_n [[T^i_n]]$, $(T^i_n)_n$ totally inaccessible times. On the other hand, since by \eqref{PDP_J=K} $J=K$, we have $D= K\cup (D\,\cap J^c)$, therefore Hypothesis \ref{H_nu}-(i) holds.

Let us now consider Hypothesis \ref{H_nu}-(ii). Taking into account  
\eqref{K_PDP}, we have to prove that  for every predictable time $S$ such that   $[[S]] \subset \{(\omega,t):X_{t-}(\omega)\in \{0,1\}\}$, 
\begin{equation}\label{PDP_muBIS} 
\nu(\{S\},\,de)=\mu(\{S\},de)\quad \textup{a.s.} 
\end{equation} 
Let $S$ be a   predictable time satisfying  $[[S]] \subset \{(\omega,t):X_{t-}(\omega)\in \{0,1\}\}$. By \eqref{mu_PDP}, $\mu(\{S\},de) = \delta_{X_S}(de)$, while from \eqref{nuPDPs} we get  $\nu(\{S\},\,de)= P(X_{S-},\,de)$. Therefore identity \eqref{PDP_muBIS} can  be rewritten as 
\begin{equation}\label{PDP_mu} 
P(X_{S-},\,de)=\delta_{X_S}(de)\quad \textup{a.s.} 
\end{equation} 
Previous identity  holds true under assumptions \eqref{cond_PDP_BIS} and \eqref{cond_PDP}, and so Hypothesis \ref{H_nu}-(ii) is established. 
 
In order to prove the validity of Hypothesis \ref{H_X_mu}, we will make use of  Lemma \ref{Ex_guida}. 
We recall that the process $X$ satisfies the stochastic differential equation  \eqref{PDP_dynamic_estended}, which gives, taking into account \eqref{nuPDPs} and Lemma \ref{L_Ban_int}, 
\begin{align}\label{X_PDP_rewr} 
X_t &= X_0 +  \int_{0}^{t}h(X_s)\,ds+ \int_{]0,\,t]}\left(\int_{\R}(e - X_{s})\,P(X_{s},\,de)\right)\lambda(X_{s})\,ds\nonumber\\ 
& 
+ \int_{]0,\,t]}(\beta(X_{s-})-X_{s-})\,dp^\ast_s+ \int_{]0,\,t]\times \R}(e-X_{s-})\,(\mu-\nu)(ds\,de). 
\end{align} 
We can show that previous equation is a particular case of  \eqref{X_SDE}. 
Indeed, we recall that, by \eqref{p_ast} and \eqref{K_PDP}, the support of the measure $d p^\ast$ is included in $K$. 
We set $B_s=s+ p^\ast(s)$ and $b(s,x)=\left(h(x)+\int_{\R}(e - x)\,\lambda(x)\,P(x,\,de)\right)\one_{K^c}(s)+(\beta(x)-x)\,\one_{K}(s)$. 
The reader can easily show that the sum of the first, second, and third integral in the right-hand side of \eqref{X_PDP_rewr} equals 
$ 
\int_0^t b(s,X_{s-})\, d B_s, 
$ 
provided we show that $\int_0^T |b(s,X_{s-})|\, d B_s$ is finite a.s. 
In fact we have 
\begin{align} \int_0^t |b(s,X_{s-})|\, d B_s 
&\leq \int_0^t |h(X_{s})|\,ds\nonumber\\ 
&+\int_{]0,\,t]}\Big |\int_{\R}(e-X_{s-})\,\lambda(X_{s-})\,P(X_{s-},\,de)\,\one_{K^c}(s)+ (\beta(X_{s-})-X_{s-})\,\one_{K}(s)\Big|\,d B_s\nonumber\\ 
&= \int_0^t |h(X_{s})|\,ds\nonumber\\ 
&+\int_{]0,\,t]}\Big |\int_{\R}(e-X_{s-})\,P(X_{s-},\,de)\,(\lambda(X_{s-})\,\one_{K^c}(s)+ \one_{K}(s))\Big|\,(ds+ d p^{\ast}(s))\nonumber\\ 
&\leq \int_0^t |h(X_{s})|\,ds+\int_{]0,\,t]}\int_{\R}|e-X_{s-}|\,\nu(ds,\,de).\label{int_2} 
\end{align} 
Recalling Lemma \ref{L_Ban_int}, and taking into account that $h$ is locally bounded, we get that 
$ 
\int_0^\cdot |b(s,X_{s-})|\, d B_s 
$ 
belongs  to  $\mathcal A_{\rm loc}^+$. 
Then, setting   $N_s=0$ and $\gamma(s,x,e)=e - x$, 
we see that $X$ is a solution to 
equation \eqref{X_SDE}. 
 
Then, by Lemma \ref{Ex_guida},  $(X, \mu)$ satisfies Hypothesis   \ref{H_X_mu} 
with decomposition $X=X^i+X^p$, where $X^i$ and $X^p$ are given respectively by \eqref{Xi_PDP} and \eqref{Xp_PDP}. 
In particular, the process $X^i$  fulfills Hypothesis \ref{H_X_mu}-1.
with $\tilde \gamma(\omega, s,e) = (e-X_{s-}(\omega))\,(1-\one_{K}(\omega, s))= (e-X_{s-}(\omega))\,\one_{\{X_{s-}(\omega)\in ]0,1[\}}(\omega,s)$. 
\endproof


\begin{lemma}\label{L_Bandini_EC0} 
We set $E = [0,\,1]$.
	Let   $(Y,U)\in 	\mathcal{L}^2 \times  \mathcal{G}^2(\mu)$ be a  solution  to the BSDE \eqref{BSDE_PDP} and  $X$
	be the piecewise deterministic Markov process   with local characteristics $(h,\lambda,P)$ satisfying \eqref{PDP_dynamic_estended}.
	Assume that $Y_t = v(t,X_t)$ for some  
	 function   $v:[0,\,T] \times \R \rightarrow \R$ such that its restriction to $[0,\,T]\times E$ is continuous. 
	Then $(X,Y)$ satisfies Hypothesis \ref{H_chain_rule_C0} with corresponding function $v$. 
\end{lemma} 
\proof	 

By Remark \ref{R_1BIS}, it follows that Hypothesis \ref{H_chain_rule_C0}-(iii) is verified.
Taking into account Remark \ref{R_2BIS},
 it follows that $Y_t=v(t, X_t)$ has a finite number of jumps, 
 so that $\sum_{s \leq T}| \Delta Y_s| < \infty$ a.s.  On the other hand, by Remark \ref{R_orth_BSDE}, $Y$ is an $(\mathcal F_t)$-orthogonal process,
so that Hypothesis \ref{H_chain_rule_C0}-(i) holds.

%

It remains to show that $v(t,X_t)$ satisfies condition \eqref{EC_C0}. 
We have 
\begin{align}\label{4.61} 
\int_{]0,\,\cdot] \times \R}|v(s,\,X_{s-}+x)-v(s,\,X_{s-})|\,\mu^X(ds\,dx) &= \sum_{0<s \leq \cdot}|v(s,\,X_{s})-v(s,\,X_{s-})|
=\sum_{s \leq \cdot}|\Delta Y_s|,
\end{align} 
by Remark \ref{R_2BIS}.
The process $Y$ takes values in the image of $[0,\,T]\times[0,1]$ with respect to $v$, which is a compact set. Therefore the jumps of $Y$ are bounded, and  \eqref{4.61} belongs to $\mathcal A_{\rm loc}^+$, see for instance the proof of Corollary at page  373 in \cite{protter}. 
\endproof 
  Finally, we apply Theorem \ref{P_ident_C0} to the present framework.
\begin{corollary}\label{C_id_PDPs} 
	Let   $(Y,U)\in 	\mathcal{L}^2 \times  \mathcal{G}^2(\mu)$ be a solution  to the BSDE \eqref{BSDE_PDP}, 
	and   $X$ the piecewise deterministic Markov process   with local characteristics $(h,\lambda,P)$ satisfying \eqref{PDP_dynamic_estended}.
	Assume that $Y_t = v(t,X_t)$ for some continuous function  $v$. 
	Assume in addition that there exists a  function $\beta: \{0,1\}\rightarrow \R$, such that 
	\begin{equation}\label{cond_PDP2} 
	P(x,\,de)\,\one_{\{x\in \{0,1\}\}}(s)=\delta_{\beta(x)}(de). 
	\end{equation} 
	Then  the random field $U$	satisfies \eqref{id_3disc} 
	with 
	\begin{align*} 
	H_s(e)&:=(U_s(e)-(v(s,e)-v(s,X_{s-}))\,\one_{\{X_{s-}\in ]0,1[\}}(s)+ U_s(e)\,\one_{\{X_{s-}\in \{0,1\}\}}(s). 
	\end{align*} 
	If in addition 
	$H_s(e)\in \mathcal G^2_{\rm loc}(\mu)$, 
	\begin{equation}\label{final_ident} 
	U_s(e)=v(s,e)-v(s,X_{s-}) \quad d\P\,\lambda(X_{s-})\,P(X_{s-},\,de)\,ds\textup{-a.e.} 
	\end{equation} 
and  there exists a predictable process $(l_s)$ such that 
\begin{equation}\label{final_ident_nud}
U_s(e)= l_s\,\one_{\{ X_{s-}\in \{0,1\}\}}(s), \quad d\P\,
\delta_{\beta(X_{s-})}(de)\,dp^\ast_s\textup{-a.e.}
\end{equation}
\end{corollary} 
\proof 
We will apply  Theorem \ref{P_ident_C0}. 
By Lemma \ref{L_Bandini_X}, 
$\mu$  satisfies Hypothesis \ref{H_nu} and $(X, \mu)$ fulfills   Hypothesis \ref{H_X_mu} with 
$\tilde \gamma(\omega,s,e)=(e-X_{s-}(\omega))\,\one_K(\omega,s)$.
Moreover, by Lemma \ref{L_Bandini_EC0}, Hypothesis \ref{H_chain_rule_C0} holds for $(X,Y)$. 
We are then in condition to apply  Theorem \ref{P_ident_C0}. 
Identity \eqref{id_3disc} holds
with 
\begin{align}\label{K_Ban} 
H_s(e)&:= U_s(e)-[v(s,X_{s-}+ \tilde \gamma(s,e))-v(s,X_{s-})]\nonumber\\ 
&=U_s(e)-[v(s,X_{s-}+ (e-X_{s-})\,
\one_K(s))-v(s,X_{s-})]\nonumber\\ 
&=[U_s(e)-(v(s,e)-v(s,X_{s-}))]\,\one_{K^c}(s)+ U_s(e)\,\one_{K}(s). 
\end{align} 
At this point we recall that 
$\nu^c(ds\,de)=\lambda(X_{s-})\, P(X_{s-},de)\,ds$, see  \eqref{nuc_PDP}. Moreover, since   $J=K$, 
\begin{equation}\label{nud_PDP} 
\nu^d(ds \,de)=\nu(ds\,de)\,\one_{K}(s)=P(X_{s-},\,de)\,dp^\ast_s=\delta_{\beta(X_{s-})}(de)\,dp^\ast_s. 
\end{equation} 
Then, since by \eqref{K_PDP} we have $K =\{(\omega,s): X_{s-}(\omega) \in \{0,1\}\}$,   \eqref{final_ident} and \eqref{final_ident_nud} are direct consequences respectively of \eqref{ID_cont_2} and \eqref{ID_disc_2}.

\endproof

\appendix 
\renewcommand\thesection{Appendix} 
\section{} 
\renewcommand\thesection{\Alph{subsection}} 
\renewcommand\thesubsection{\Alph{subsection}}

\subsection{Technical results related to the  hypothesis on the underlying process $X$}\label{Appendix_A}
The results below are related to Hypothesis \ref{H_X_mu} concerning $(X,\mu)$, with $X$ being a c\`adl\`ag process  and $\mu$ an integer-valued random measure; they will be extensively used in Appendix \ref{Appendix B}.
\begin{remark}\label{R_exhseq} 
\begin{itemize}
\item[(i)]
	Given a predictable thin set $A$,
	there exists a 
	 sequence of predictable times $(R_n)_n$  with disjoint graphs, such that $A = \cup_n [[R_n]]$, up to an evanescent set,
		see 
	Proposition 2.23, Chapter I, in \cite{JacodBook}.
	\item [(ii)]
	Since $\{\Delta X^p \neq 0\}$ is a predictable thin set (see the comments after Definition 7.39 in \cite{chineseBook}),
		by item (i) there exists a 
	 sequence of predictable times 
	 exhausting the jumps of $X^p$, up to an evanescent set.
	\end{itemize}
\end{remark} 
\begin{proposition}\label{R_NEW} 
	Let $X$ be a c\`adl\`ag adapted  process  
	with decomposition 
	$X=X^i+X^p$, where 
	$X^i$ (resp. $X^p$) is a c\`adl\`ag 
	quasi-left continuous 
	adapted process (resp. c\`adl\`ag  predictable  process). Then the two properties below hold. 
	\begin{itemize} 
		\item[(i)] 
		$\Delta X^p\,\one_{\{\Delta X^i \neq 0\}}=0$ and $\Delta X^i\,\one_{\{\Delta X^p \neq 0\}}=0$, up to an evanescent set. 
		 
		\item[(ii)]$\{\Delta X \neq 0\}$ is the disjoint union (up to an evanescent set) of the random sets $\{\Delta X^p \neq 0\}$ and  $\{\Delta X^i \neq 0\}$.
	\end{itemize} 
\end{proposition}	 
\proof	 
(i)	
By Remark \ref{R_exhseq}-(ii), there exist a sequence of predictable times $(T_n^p)_n$  that exhausts the jumps of $X^p$, up to an evanescent set. 
Moreover, recalling  
 Proposition 2.26, Chapter I, in \cite{JacodBook}, there exist a sequence of 
  totally inaccessible times $(T_n^i)_n$) that exhausts the jumps of 
 $X^i$. 
On the other hand, $\Delta X^p_{T^i_n} = 0$ a.s. for every $n$, see 
Proposition 2.24, Chapter I, in \cite{JacodBook}
  (resp.  $\Delta X^i_{T^p_n} =0$ a.s. for every $n$, see 
  Definition 2.25, Chapter I, in \cite{JacodBook}), so that, up to an evanescent set,  
\begin{align*} 
\Delta X^i\,\one_{\{\Delta X^p \neq 0\}}=\Delta X^i\,\one_{\cup_n [[T_n^p]]}=0,\qquad 
\Delta X^p\,\one_{\{\Delta X^i \neq 0\}}=\Delta X^p\,\one_{\cup_n [[T_n^i]]}=0. 
\end{align*} 
	
(ii) We have, again up to an evanescent set, 
\begin{eqnarray*}
&&\{\Delta X\neq 0 \}=\{(\Delta X^i + \Delta X^p) \neq 0 \}\\ 
&&=(\{(\Delta X^i + \Delta X^p) \neq 0 \} \,\cap \,\{\Delta X^p =0 \})\cup\, (\{(\Delta X^i + \Delta X^p) \neq 0 \} \,\cap \,\{\Delta X^p \neq 0 \})\\
&&=(\{\Delta X^i  \neq 0 \} \,\cap \,\{\Delta X^p =0 \})\cup\, \{\Delta X^p \neq 0 \}\\
&&=\{\Delta X^i\neq 0\} \cup \{\Delta X^p\neq 0 \},
\end{eqnarray*}
where the third equality follows from the second statement in  item (i) of the Proposition.
We observe that the intersection of $\{\Delta X^i \neq 0\}$ and $\{\Delta X^i \neq 0\}$ is evanescent because of   item (i).

\endproof 

\begin{proposition}\label{P_1} 
	Let $X$ be a c\`adl\`ag adapted process with decomposition 
	$X=X^i+X^p$, where 
	$X^i$ (resp. $X^p$) is a c\`adl\`ag 
	quasi-left continuous 
	adapted process (resp. c\`adl\`ag  predictable  process). 
	Then 
		$$\{(\omega,t) : \nu^X(\omega, \{t\}\times \R) >0\}=\{\Delta X^p \neq 0\}.
		$$
\end{proposition} 
\proof 
$\{\Delta X \neq 0\}$ is the support of the random measure $\mu^X$ (see e.g. 
Proposition 1.16, Chapter II,  in \cite{JacodBook}). 
By 
Theorem 11.14 in \cite{chineseBook}, 
the  predictable support of $\{\Delta X \neq 0\}$ is given by $\{(\omega, t) : \nu^X(\{t\}\times \R) >0\}$.	 
 
It remains to prove that the predictable support of $\{\Delta X\neq 0\}$ equals  $\{\Delta X^p\neq 0\}$. 
By Proposition \ref{R_NEW}-(ii), $\{\Delta X \neq 0\}$ is the disjoint union (up to an evanescent set) of  $\{\Delta X^p \neq 0\}$ and  $\{\Delta X^i \neq 0\}$. 
Since $X^i$ is a c\`adl\`ag quasi-left continuous process, by 
Proposition 2.35, Chapter I, in \cite{JacodBook}, we know that  the  predictable support of $\{\Delta X^i \neq 0\}$ is evanescent. 
Then, by the definition of predictable support, see 
Definition 2.32, Chapter I, in \cite{JacodBook}, taking into account the additivity of the predictable projection operator, we have 
${}^{p}\left(\one_{\{\Delta X \neq 0\}}\right) = \one_{\{\Delta X^p \neq 0\}}$,
and this concludes the proof. 
\endproof 
\begin{proposition}\label{P_Tnk} 
		Let $X$ be a c\`adl\`ag adapted process with decomposition 
	$X=X^i+X^p$, where 
	$X^i$ (resp. $X^p$) is a c\`adl\`ag 
	quasi-left-continuous 
	adapted process (resp. c\`adl\`ag  predictable  process). 
	Let $(S_n)_n$ be a sequence of predictable times exhausting the jumps of $X^p$. 
	Then 
	\begin{equation} 
	\nu^X(\{S_{n}\},dx)= \mu^X(\{S_{n}\},dx)\,\,\textup{for any}\,\, n,\,\, \textup{a.s.} 
	\end{equation} 
\end{proposition} 
\proof 
Let us fix $n \in \N$.
We need to show the existence of a $\P$-null set $\mathcal N$ such that, for every $\omega \notin \mathcal N$, we have  
\begin{equation}\label{final_eq}
\int_\R \one_{E_m}(x)\,\nu^X(\{S_{n}\}, dx)=	\int_\R \one_{E}(x)\,\mu^X(\{S_{n}\}, dx)
\end{equation}
for every real Borel set $E$.
Let  $(E_m)_m$ be a  sequence of measurable subsets of $\R$ which is 
a $\pi$-class generating $\mathcal{B}(\R)$. 
Since $X^i$ is  a c\`adl\`ag quasi-left-continuous adapted process and $S_n$ is a predictable time, then $\Delta X^i_{S_n} =0$ a.s., see again
Definition 2.25, Chapter I, in \cite{JacodBook}. This implies that $\Delta X_{S_{n}}=\Delta X^p_{S_{n}}$ a.s. by Hypothesis \ref{H_X_mu}-2. 
Consequently, for every $m$ we have 
\begin{equation}\label{E_1} 
\one_{E_m}(\Delta X^p_{S_{n}}) = \one_{E_m}(\Delta X_{S_{n}})= \int_\R \one_{E_m}(x)\,\mu^X(\{S_{n}\}, dx)\,\,\textup{a.s.} 
\end{equation} 
On the other hand, by point (b) of Proposition 1.17, Chapter II, in \cite{JacodBook}
 and \eqref{E_1} we have 
\begin{align*} 
\int_\R \one_{E_m}(x)\,\nu^X(\{S_{n}\}, dx) &= \sper{\int_\R  \one_{E_m}(x)\,\mu^X(\{S_{n}\}, dx)\Big|\mathcal{F}_{S_{n} -}}\\ 
&= \sper{\one_{E_m}(\Delta X^p_{S_{n}})\Big|\mathcal{F}_{{S_{n}} -}}\\ 
&=\one_{E_m}(\Delta X^p_{S_{n}})\,\,\textup{a.s.}
\end{align*} 
By \eqref{E_1}, there exists a 
null set $\mathcal{N}_m$ such that 
$$ 
\int_\R \one_{E_m}(x)\,\nu^X(\{S_{n}\}, dx)=	\int_\R \one_{E_m}(x)\,\mu^X(\{S_{n}\}, dx)\,\,	\textup{for every} \,\,\omega \notin \mathcal{N}_m. 
$$ 
Define $\mathcal{N} = \cup_m \mathcal{N}_m$, then 
$$ 
\int_\R \one_{E_m}(x)\,\nu^X(\{S_{n}\}, dx)=	\int_\R \one_{E_m}(x)\,\mu^X(\{S_{n}\}, dx)\,\,	\textup{for every}\,\,m \,\,\textup{and} \,\,\omega \notin \mathcal{N}. 
$$ 
Then \eqref{final_eq} follows by a monotone class argument, see 
Theorem 1.1, Chapter 1, \cite{kallenberg}. 
\endproof	


 \begin{proposition}\label{P_M_mu} 
	Let  $Y$ be 
	a c\`adl\`ag  adapted  process such that $(Y, \mu)$ satisfies 
	Hypothesis 
\ref{H_X_mu}-1.
	Then, 	 there exists a null set $\mathcal N$ such that,  for every   Borel  function $\varphi: [0,\,T] \times \R \rightarrow \R_+$ satisfying  $\varphi(s,0)=0$ for every $s \in [0,\,T]$, we have 
	\begin{equation}\label{E_H_(ii)} 
	\sum_{0 <s \leq T} 
	\varphi(s,\Delta Y_s(\omega)) = \int_{]0,\,T]\times \R} 
	\varphi(s,\tilde \gamma(\omega,s, e))\,\mu(\omega,ds\,de), 
	\quad \omega \notin \mathcal N. 
	\end{equation} 
\end{proposition} 
\proof 
Taking into account that $\{\Delta Y\neq 0\}\subset D$ and the fact that $\varphi(s,0) =0$, it will be enough to prove that 
\begin{equation}\label{E_H_(ii)BIS} 
\sum_{0<s \leq T} 
\varphi(s,\Delta Y_s(\omega))\,\one_{D}(\omega,s) = \int_{]0,\,T]\times \R} 
\varphi(s,\tilde \gamma(\omega,s, e))\,\mu(\omega,ds\,de), 
\quad \omega \notin \mathcal N, 
\end{equation} 
for every   Borel function $\varphi: [0,\,T] \times \R \rightarrow \R_+$. 
 
Let  $(I_m)_m$ be a sequence of subsets of $[0,\,T]\times \R$, which is a $\pi$-system generating $\mathcal B([0,\,T]) \otimes \mathcal{B}(\R)$. 
Setting $\varphi_m(s,x)=\one_{I_m}(s,x)$, for every $m$, we will show that 
\begin{equation}
\label{id_as} 
\sum_{0<s \leq T} 
\varphi_m(s,\Delta Y_s)\, \one_{D}(\cdot, s)
= \int_{]0,\,T]\times \R} 
\varphi_m(s,\tilde \gamma(\cdot,s,e))\, 
\mu(\cdot,ds\,de),\quad \textup{a.s.} 
\end{equation} 
Let $n \in \N$ be fixed. 
In order to establish \eqref{id_as},  it is enough to prove 
\begin{equation} 
\label{id_asBIS} 
\sum_{0<s \leq T} 
\varphi_m(s,\Delta Y_s)\,\,\one_{\tilde{\Omega}_n}(\cdot,s)\, \one_{D}(\cdot, s)
= \int_{]0,\,T]\times \R} 
\varphi_m(s,\tilde \gamma(\cdot,s,e))\, 
\,\one_{\tilde{\Omega}_n}(\cdot, s)\,\mu(\cdot,ds\,de),\quad \textup{a.s.,} 
\end{equation} 
where $\tilde{\Omega}_n$ are the sets introduced in Definition \ref{D_DoleansMeas}.
Let us consider a bounded, $\mathcal{F}$-measurable function $\phi: \Omega \rightarrow \R_{+}$. Identity \eqref{id_asBIS} holds if 
we  show that the expectations of both sides 
against $\phi$ are equal. 
Using  Hypothesis 
\ref{H_X_mu}-1, we write 
\begin{align*} 
&\E\bigg[\phi\,\int_{]0,\,T]\times \R} 
\varphi_m(s,\tilde \gamma(\cdot,s,e))\,\one_{\tilde{\Omega}_n}(\cdot, s)\,\mu(\cdot,ds\,de)\bigg]\\ 
&=\int_{\Omega \times ]0,\,T]\times \R} dM^\P_{\mu}(\omega,s,e)\,\phi(\omega)\,\varphi_m(s,\tilde\gamma(\omega,s,e))\,\one_{\tilde{\Omega}_n}(\omega, s)\,\\ 
&= \int_{\Omega \times ]0,\,T]} d M^\P_{\mu}(\omega,s, y)\,\phi(\omega)\,\varphi_m(s,\Delta Y_s(\omega))\,\one_{\tilde{\Omega}_n}(\omega, s)\\ 
&=\E\bigg[\phi\,\int_{]0,\,T]\times \R} 
\varphi_m(s,\Delta Y_s (\cdot))\,\one_{\tilde{\Omega}_n}(\cdot, s)\,\mu(\cdot,ds\,dy)\bigg]\\ 
&=\E\bigg[\phi\,\int_{]0,\,T]\times \R}\sum_{s> 0}
\varphi_m(s,\Delta Y_s (\cdot))\,\one_{\tilde{\Omega}_n}(\cdot, s)\,\one_{D}(\cdot,s)\,\delta_{(s,\beta_s(\cdot))}(dt\,dx)\bigg]\\ 
&= \E\Bigg[\phi\sum_{0<s \leq T}
\one_{D}(\cdot,s)\,\varphi_m(s,\Delta Y_s(\cdot))\,\one_{\tilde{\Omega}_n}(\cdot, s)\Bigg], 
\end{align*} 
where we have used the form of $\mu$ given  in Proposition 1.14, Chapter II,  in \cite{JacodBook}, i.e.
\begin{equation}\label{mubeta}
\mu(dt\,dy) = \sum_{s \geq 0}\one_{D}(s,\omega)\,\delta_{(s,\beta_s(\omega))}(dt\,dy).
\end{equation}
Therefore, there exists a $\P$-null set $\mathcal N_m$ such that 
$$ 
\sum_{0<s \leq T} 
\varphi_m(s,\Delta Y_s(\omega))\,\one_{D}(\omega,s) \,\one_{\tilde{\Omega}_n}(\omega,s)
= \int_{]0,\,T]\times \R} 
\varphi_m(s,\tilde \gamma(\omega, s, e))\,\one_{\tilde{\Omega}_n}(\omega,s)
\mu(\omega,ds\,de),\quad \omega \notin \mathcal{N}_m. 
$$ 
Define $\mathcal{N} = \cup_m \mathcal{N}_m$, then for $\varphi = \varphi_m$ for every $m$ we have 
$$ 
\sum_{0<s \leq T} 
\varphi(s,\Delta Y_s(\omega))\,\one_{D}(\omega,s) \,\one_{\tilde{\Omega}_n}(\omega,s)
= \int_{]0,\,T]\times \R} 
\varphi(s,\tilde \gamma(\omega, s, e))\,\one_{\tilde{\Omega}_n}(\omega,s) 
\mu(\omega,ds\,de),\quad \omega \notin \mathcal{N}. 
$$ 
By a monotone class argument (see Theorem 2.3 in   \cite{jacodprotter}) the identity holds for every measurable bounded function  $\varphi: [0,\,T] \times \R \rightarrow \R$, and therefore, using monotone convergence theorem,  for every positive measurable function $\varphi$ on $[0,\,T] \times \R$ as well. 
\endproof

\subsection{Proofs of the technical results stated in Section \ref{Sec_2.2}}\label{Appendix B}
 \subsubsection{Proof of Proposition \ref{P_identiy_measures}}\label{Sec_Proof_P_identiy_measures}
Let $\varphi:  [0,\,T] \times \R \rightarrow \R_+$. 
Taking into account Proposition \ref{R_NEW}-(i) and the fact that $\varphi(s,0)=0$, we have, for almost all $\omega$, 
\begin{align*} 
\sum_{0<s \leq T} \varphi(s,\Delta X_s(\omega)) 
&=\sum_{0<s \leq T} \varphi(s,\Delta X^i_s(\omega)+\Delta X^p_s(\omega))\,\one_{\{\Delta X^p = 0\}}(\omega,s) \\ 
&+ \sum_{0<s \leq T} \varphi(s,\Delta X^i_s(\omega)+\Delta X^p_s(\omega))\,\one_{\{\Delta X^p \neq 0\}}(\omega,s)\\ 
&=\sum_{0<s \leq T} \varphi(s,\Delta X^i_s(\omega))\,\one_{\{\Delta X^p = 0\}}(\omega,s) + \sum_{0<s \leq T} \varphi(s,\Delta X^p_s(\omega))\,\one_{\{\Delta X^p \neq 0\}}(\omega,s)\\ 
&=\sum_{0<s \leq T} \varphi(s,\Delta X^i_s(\omega)) + \sum_{0<s \leq T} \varphi(s,\Delta X^p_s(\omega)). 
\end{align*} 
By Proposition \ref{P_M_mu} applied to $Y= X^i$,  
there exists a null set $\mathcal N$ such that, for every $\omega \notin \mathcal N$, previous expression gives 
\begin{align*} 
\int_{]0,\,T]\times \R}\varphi(s,x)\,\mu^X(\omega,ds\,dx)=\int_{]0,\,T]\times \R}\varphi(s,\tilde  \gamma(\omega,s,e))\,\mu(\omega,ds\,de)+ \sum_{0<s \leq T} \varphi(s,\Delta X^p_s(\omega)). 
\end{align*} 
The second part of the statement holds decomposing  $\varphi = \varphi^+ - \varphi^-$. 
\qed

\subsubsection{Proof of Proposition \ref{L_ident_mu_muX}}\label{Sec_Proof_L_ident_mu_muX}
Clearly the result holds if we show that $\varphi$ verifies \eqref{Id_intstoch} under  one of the two following assumptions: 
\begin{itemize} 
	\item[(i)] $|\varphi| \star \mu^X \in \mathcal A_{\rm loc}^+$, 
	\item[(ii)] $|\varphi|^2 \star \mu^X \in \mathcal A_{\rm loc}^+$. 
\end{itemize}	 
By localization arguments, it is enough to show it when $|\varphi| \star \mu^X \in \mathcal A^+$, $|\varphi|^2 \star \mu^X \in \mathcal A^+$. 
Below we will consider the first case, the second case will follow from the first one by approaching in $\mathcal L^2(\mu^X)$ the function $\varphi$ with $\varphi_{\varepsilon}(s,x):=\varphi(s,x)\,\one_{\varepsilon < |x| \leq 1/\varepsilon}\,\one_{s \in [0,\,T]}$. Indeed, $\varphi_\varepsilon(s,x) \star \mu \in \mathcal A^+$, by Cauchy-Schwarz inequality,
taking into account the fact that $\mu^X$, restricted to $\varepsilon \leq |x| \leq 1/\varepsilon$, is finite, since $\mu^X$ is $\sigma$-finite on $[0,\infty) \times \R$.  
 
Let us define 
\begin{align} 
M_t &:= \int_{]0,\,t]\times \R}\varphi(\cdot, s,x)\,(\mu^X-\nu^X)(ds\,dx),\nonumber\\ 
N_t &:= \int_{]0,\,t]\times \R}\varphi(\cdot, s,\tilde \gamma(\cdot,s,e)) \,(\mu-\nu)(ds\,de).\label{N} 
\end{align} 
Notice that  the processes $M$ and $N$ are purely discontinuous  local martingales, see e.g. 
Definition 1.27, point b), Chapter II, in \cite{JacodBook}. 
We have to prove that $M$ and $N$ are indistinguishable. To this end, by 
Corollary 4.19, Chapter I,  in \cite{JacodBook}, it is enough to prove that $\Delta M = \Delta N$, up to an evanescent set. 
Observe that 
\begin{align} 
\Delta M_s &= \int_{\R}\varphi(\cdot,s,x)\,(\mu^X-\nu^X)(\{s\},\,dx)\nonumber\\ 
&=\int_{\R}\varphi(\cdot, s,x)\,(1-\one_{J}(\cdot, s))\, (\mu^X-\nu^X)(\{s\},\,dx)\nonumber\\ 
&+ \int_{\R}\varphi(\cdot, s,x)\,\one_{J}(\cdot, s)\, (\mu^X-\nu^X)(\{s\},\,dx),\label{DeltaM} 
\end{align} 
and \begin{align} 
\Delta N_s &= \int_{\R}\varphi(\cdot, s,\tilde \gamma(\cdot, s,e))\,(\mu-\nu)(\{s\},\,de)\nonumber\\ 
&= \int_{\R}\varphi(\cdot, s,\tilde \gamma(\cdot, s,e))\,\one_{J}(\cdot, s)\,(\mu-\nu)(\{s\},\,de)\nonumber\\ 
&+ \int_{\R}\varphi(\cdot, s,\tilde \gamma(\cdot, s,e))\,(1-\one_{J}(\cdot,s))\,(\mu-\nu)(\{s\},\,de).\label{DeltaN} 
\end{align} 
By 	definition of  $J$, for every $\omega$ and every $s$ we have 
\begin{equation}\label{nu_1-J} 
\nu(\omega,\{s\},de)\,(1-\one_{J}(\omega, s)) =0. 
\end{equation} 
Moreover, since $J$ is a predictable thin set (see Remark \ref{R_pred_supp}), there exists 
a sequence of predictable times $(R_n)_n$  with disjoint graphs, such that $J = \cup_n [[R_n]]$, see Remark \ref{R_exhseq}-(i). We recall that Hypothesis \ref{H_nu}-(i) implies that $J=K$ up to an evanescent set, see Remark \ref{R_Hp_3.1}. 
By this fact, and taking into account  Hypothesis \ref{H_nu}-(ii),  there exists a null set $\mathcal N$, such that, 
for every $n \in \N$, $\omega \notin \mathcal N$, 
$$ 
\mu(\omega,\{R_n(\omega)\},de)\,\one_{J}(\omega, s) =\nu(\omega,\{R_n(\omega)\},de)\,\one_{J}(\omega, s). 
$$ 
By additivity, it follows that for every $\omega \notin \mathcal N$, for every $s \in [0,\,T]$, 
\begin{equation} 
\label{munuJ} 
\mu(\omega,\{s\},de)\,\one_{J}(\omega, s) =\nu(\omega,\{s\},de)\,\one_{J}(\omega, s). 
\end{equation} 
 
On the other hand,  $\{\Delta X^p \neq 0\} \subset J$ by Hypothesis \ref{H_X_mu}-2. 
Recalling that $\{\Delta X^p \neq 0\} 
=\{(\omega,s): \nu^X(\{s\}\times \R)>0\}$ 
(see Proposition \ref{P_1}), 
for almost every $\omega$, for every $s \in [0,\,T]$, we have 
\begin{equation}\label{rightside} 
\nu^X(\omega, \{s\},dx)\,\one_{J}(\omega, s)=\nu^X(\omega, \{s\},dx)\,\one_{\{\Delta X^p \neq 0 \}}(\omega, s), 
\end{equation} 
so that 
\begin{equation}\label{nuX_1-J} 
\nu^X(\omega, \{s\},dx)\,(1-\one_{J}(\omega, s))= 
\nu^X(\omega, \{s\},dx)\,(1-\one_{\{\Delta X^p\neq 0\}}(\omega, s)) 
=0. 
\end{equation} 
Now notice that there always exists a sequence of predictable 
times  exhausting the jumps of $X^p$, up to an evanescent set, see Remark \ref{R_exhseq}-(ii). 
By means of Proposition \ref{P_Tnk} 
we can prove, similarly as we did in order to establish \eqref{munuJ},   that for every $\omega \notin \mathcal N$ ($\mathcal N$ possibly enlarged), for every $s \in [0,\,T]$, 
\begin{align}\label{third} 
&\mu^X(\omega, \{s\},dx)\,\one_{\{\Delta X^p \neq 0\}}(\omega, s)=\nu^X(\omega, \{s\},dx)\,\one_{\{\Delta X^p \neq 0\}}(\omega, s). 
\end{align} 
Finally, we notice that 
\begin{equation}\label{E_n}
\mu^X(\omega, \{s\},dx)\,\one_{J}(\omega, s)=\mu^X(\omega, \{s\},dx)\,\one_{J \cap \{\Delta X \neq 0\}}(\omega, s). 
\end{equation}
Taking into account that   $X^i$ is a  c\`adl\`ag  quasi-left-continuous process, by 
Definition 2.25, Chapter I, in \cite{JacodBook} we have $\Delta X^i_{R_n}=0$ for every $n$, so that
\begin{align*} 
J \cap \{\Delta X \neq 0\}& = (\cup_n [[R_n]] \cap \{\Delta X^i \neq 0\}) \cup (\cup_n [[R_n]] \cap \{\Delta X^p \neq 0\})\\ 
& =  \cup_n [[R_n]] \cap \{\Delta X^p \neq 0\}=  \{\Delta X^p \neq 0\}. 
\end{align*} 
This, together with \eqref{E_n}, implies that,  for every $\omega \notin \mathcal N$   and for every $s \in [0,\,T]$, 
\begin{align}\label{leftside} 
\mu^X(\omega, \{s\},dx)\,\one_{J}(\omega, s)&=\mu^X(\omega, \{s\},dx)\,\one_{J \cap \{\Delta X \neq 0\}}(\omega, s)=\mu^X(\omega, \{s\},dx)\,\one_{\{\Delta X^p \neq 0\}}(\omega, s). 
\end{align} 
Collecting \eqref{rightside}, \eqref{third} and \eqref{leftside} we conclude that for every $\omega \notin \mathcal N$, for every $s \in [0,\,T]$, 
\begin{equation}\label{muXnuX_J} 
\mu^X(\omega, \{s\},dx)\,\one_{J}(\omega, s)=\nu^X(\omega, \{s\},dx)\,\one_{J}(\omega, s). 
\end{equation} 
Therefore, for every $\omega \notin \mathcal N$, for every $s \in [0,\,T]$, taking into account \eqref{nu_1-J}, \eqref{munuJ},  \eqref{nuX_1-J}, \eqref{muXnuX_J},  expressions  \eqref{DeltaM} and \eqref{DeltaN} become 
\begin{align} 
\Delta M_s &= \int_{\R}\varphi(\cdot, s,x)\,(1-\one_{J}(\cdot,s))\,\mu^X(\{s\},\,dx),\label{DeltaM_final}\\ 
\Delta N_s &= \int_{\R}\varphi(\cdot, s,\tilde \gamma(\cdot, s,e))\,(1-\one_{J}(\cdot, s))\,\mu(\{s\},\,de).\label{DeltaN_final} 
\end{align} 
 \\
Now let us prove that, for every $s \in [0,\,T]$,   $\Delta M_s(\omega)= \Delta N_s(\omega)$ for every $\omega \notin \mathcal N$, namely up to an evanescent set. 
Set 
$$ 
\varphi_s(\omega, t,x):= \varphi(\omega, t,x)\,\,(1-\one_{J}(\omega,t))\,\one_{\{s\}}(t), 
$$ 
then $\Delta M_s$ and $\Delta N_s$ can be rewritten as 
\begin{align*} 
\Delta M_s(\omega) &= \int_{[0,\,T]\times \R}\varphi_s(\omega,t,x)\,\mu^X(\omega, dt\,dx),\\ 
\Delta N_s(\omega) &= \int_{[0,\,T]\times \R}\varphi_s(\omega, t,\tilde \gamma(\omega,t,e))\,\mu(\omega, dt\,de). 
\end{align*} 
Then, Proposition \ref{P_identiy_measures} and Remark \ref{R_phi_real_val} applied  to the process $\varphi_s$ 
implies that (possibly enlarging the null set $\mathcal N$), 
\begin{equation*} 
\int_{]0,\,T]\times \R} \varphi_s(\omega, t,x)\,\mu^X(\omega,dt\,dx)= \int_{]0,\,T]\times \R} \varphi_s(t,\tilde \gamma(\omega,t,e))\,\mu(\omega,dt\,de) + V^{\varphi_s}(\omega) 
\end{equation*} 
for every $\omega \notin \mathcal N$,  or, equivalently, that 
\begin{equation*} 
\int_{\R}\varphi(\omega,s,x)\,\mu^X(\omega,\{s\},\,dx)= \int_{\R}\varphi(\omega,s,\tilde \gamma(\omega,s,e))\,\mu(\omega,\{s\},\,de) + V^{\varphi_s}(\omega), 
\end{equation*} 
for every $\omega \notin \mathcal N$, where 
\begin{equation}\label{Vvarphis} 
V^{\varphi_s}(\omega) 
= \sum_{t \leq T}\varphi_s(\omega, t,\Delta X^p_t(\omega)) 
=\varphi(\omega, s,\Delta X^p_{s}(\omega))\,\,\one_{J^c \cap \{\Delta X^p \neq 0\} }(\omega,s). 
\end{equation} 
Recalling that $\{\Delta X^p \neq 0 \}\subset J$ by Hypothesis \ref{H_X_mu}-2., 
 it straightly follows from \eqref{Vvarphis} that $V^{\varphi_s}(\omega)$ is zero. 
In particular,  up to an evanescent set, we have 
\begin{equation*} 
\int_\R\varphi(\omega, s,x)\,\mu^X(\omega, \{s\},\,dx)= \int_\R\varphi(s,\tilde \gamma(\omega, s,e))\,\mu(\omega, \{s\},\,de), 
\end{equation*} 
in other words $\Delta M=\Delta N$ up to an evanescent set, and this concludes the proof. 
\qed

\subsubsection{Proof of Lemma \ref{Ex_guida}}\label{Sec_Proof_Ex_guida}
Since $N$ is continuous, it straight follows from \eqref{Xp} that 
\begin{equation}\label{DeltaXp} 
\Delta X^p_s = b(s,X_{s-})\,\Delta B_s. 
\end{equation} 
We remark that $X^i$ in \eqref{Xi} has the same expression as $N$ defined in \eqref{N} where the integrand $\varphi(\omega,s, \tilde \gamma(\omega, s, e))$ is replaced by $\gamma(s,X_{s-}(\omega),e)$. 
We recall that  Hypothesis \ref{H_nu}-(i) implies that $J=K$ up to an evanescent set, see   Remark \ref{R_Hp_3.1}. 
Similarly as for 
\eqref{DeltaN_final}, 
we get 
\begin{equation} 
\Delta X^i_s = \int_{\R} \gamma(s,X_{s-},e)\,\,(1-\one_{K}(s))\,\mu(\{s\},\,de).\label{gammatilde} 
\end{equation} 
Since by Hypothesis \ref{H_nu}, up to an evanescent set,  
$D\setminus K = \cup_n [[T_n^i]]$, $(T_n^i)_n$ being a sequence of totally inaccessible times with disjoint graphs, recalling \eqref{mubeta}, 
\eqref{gammatilde} can be rewritten as 
\begin{align} 
\Delta X^i_s(\omega) =  \gamma(s,X_{s-}(\omega),\beta_s(\omega))\,\,\one_{\cup_n [[T_n^i]]}(\omega, s).\label{DeltaXi} 
\end{align} 
 
We can easily show that 
$X^p$ and $X^i$ are respectively  a c\`adl\`ag   predictable process and a c\`adl\`ag  quasi-left-continuous  adapted process. 	 
The fact that  $X^p$ is predictable straight follow from \eqref{Xp}. Concerning $X^i$,  let $S$ be a predictable time; it is enough to prove that  $\Delta X^i_S\,\one_{\{S< \infty\}} = 0$ a.s., see 
Definition 2.25, Chapter I, in \cite{JacodBook}. Identity  \eqref{DeltaXi} gives 
\begin{align}\label{DeltaXi_D} 
\Delta X^i_S(\omega)\,\one_{\{S< \infty\}} =  \gamma(S,X_{S-}(\omega),\beta_S(\omega))\,\,\one_{\cup_n [[T_n^i]]}(\omega, S(\omega))\,\one_{\{S< \infty\}}. 
\end{align} 
Since  the   graphs of the totally inaccessible times $T_n^i$ are disjoint, $\one_{\cup_n [[T_n^i]]}(\omega, S(\omega))\,\one_{\{S< \infty\}}= \sum_n \one_{[[T_n^i]]}(\omega, S(\omega))\,\one_{\{S< \infty\}}$, and  the conclusion follows by  
 Remark \ref{R_H_mu}-(iii), 
 taking into account that $S$ is a predictable time.
 
The process $X^p$ in \eqref{Xp} satisfies Hypothesis \ref{H_X_mu}-2. 
Indeed, by \eqref{DeltaXp} we have 
\begin{equation}\label{DeltaXp_D} 
\{\Delta X^p \neq 0\}\subset \{\Delta B \neq 0\}\subset J.
\end{equation} 
 
Finally, we show that the process $X^i$ in \eqref{Xi} fulfills Hypothesis \ref{H_X_mu}-1. 
with   $\tilde \gamma(\omega, s,e) = \gamma(s,X_{s-}(\omega),e)\,(1-\one_{K}(\omega, s))$. 
First, the fact that $\{\Delta X^i\neq 0\}\subset D$ directly follows from \eqref{gammatilde}. 
To prove $\Delta X^i_s(\omega) = \tilde \gamma(\omega, s,\cdot)$, $dM^{\P}_{\mu}(\omega,s)$-a.e. it is enough to show that 
$$ 
\sper{\int_{]0,\,T]\times \R}\mu(\omega, ds\,de)\,|\tilde \gamma(\omega, s,e)-\Delta X^i_s(\omega)|\,\one_{\tilde \Omega_n}(\omega,s)}=0. 
$$ 
To establish this, we see that by the structure of $\mu$ it follows that, for every $n \in \N$, 
\begin{align*} 
\sper{\int_{]0,\,T]\times \R}\mu(\omega, ds\,de)\,|\tilde \gamma(\omega, s,e)-\Delta X^i_s(\omega)|\one_{\tilde \Omega_n}(\omega,s)} \leqslant
\sum_{s \in ]0,\,T]} \sper{\one_{D}(\cdot,s)\,|\tilde \gamma(\cdot, s,\beta_s(\cdot))-\Delta X^i_s(\cdot)|}, 
\end{align*} 
which vanishes taking into account \eqref{DeltaXi}. 
\qed

\makeatletter 
\@addtoreset{equation}{subsection} 
\def\theequation{\thesection.\arabic{equation}} 
\makeatother 
 
\theoremstyle{plain} 
\newtheorem{Theorem}{Theorem}[subsection] 
\newtheorem{Lemma}[Theorem]{Lemma} 
\newtheorem{Proposition}[Theorem]{Proposition} 
\newtheorem{Corollary}[Theorem]{Corollary} 
 
\newtheorem{Definition}[Theorem]{Definition} 
\newtheorem{Hypothesis}[Theorem]{Hypothesis} 
\theoremstyle{remark} 
\newtheorem{Remark}[Theorem]{Remark} 
\newtheorem{Example}[Theorem]{Example}

\small 
 
\paragraph{Acknowledgements.} 
The authors wish to thank the Referees for their valuable comments 
and criticisms, which have motivated us to strongly  improve the first submitted version  of
the paper.
The first named author benefited from the support of the Italian MIUR-PRIN 2010-11 ``Evolution differential problems: deterministic and stochastic approaches and their interactions''.
The second named author  benefited 
partially from the support of the ``FMJH Program Gaspard Monge in optimization 
and operation research'' (Project 2014-1607H). 

\small 
 
\addcontentsline{toc}{chapter}{Bibliography} 
\bibliographystyle{plain} 


\bibliography{BiblioLivreFRPV_TESI} 

\end{document}